\documentclass[12pt, a4paper, twoside]{amsart}
\usepackage{amsmath}
\usepackage{amssymb}
\usepackage{amsfonts}
\usepackage{a4}
\usepackage[latin1]{inputenc}

\renewcommand{\baselinestretch}{1.2}

\numberwithin{equation}{section}

\newtheorem{theorem}{Theorem}[section]
\newtheorem{lemma}[theorem]{Lemma}

\theoremstyle{definition} 

\theoremstyle{plain} 

\newtheorem{proposition}[theorem]{Proposition}

\newtheorem{fact}[theorem]{Fact}
\newtheorem{conclusion}[theorem]{Conclusion}

\newtheorem{claim}[theorem]{Claim}

\newtheorem*{maintheorem*}{Main Theorem} 
\newtheorem*{conjecture*}{Conjecture} 
\newtheorem{definition}[theorem]{Definition}
\newtheorem{choice}[theorem]{Choice}
\newtheorem{hypothesis}[theorem]{Hypothesis}
\newtheorem{remark}[theorem]{Remark}

\theoremstyle{remark}  

\newcommand{\nc}{\newcommand}
\nc{\nothing}[1]{}
\nc{\added}[1]{#1} 
\nc{\comment}[1]{#1} 

\nothing{   
\nc{\dom}{{\rm dom}}
\nc{\card}{{\rm card}}
\nc{\lh}{{\rm lh}}
\nc{\lgg}{{\rm lg}}
\nc{\rge}{\mbox{\rm range}}
\nc{\cf}{{\rm cf}}
\nc{\nex}{\mbox{\rm next}}
\nc{\uhr}{\restriction}
\nc{\supt}{{\rm supt}}
\nc{\supp}{{\rm supp}}
\nc{\Lim}{{\rm Lim}}
\nc{\Leb}{{\rm Leb}}
\nc{\modd}{{\rm mod}}
\nc{\RO}{{\rm RO}}
\nc{\prob}{{\rm Prob}}
}

\nc{\On}{{\rm On}}

\nc{\nco}{\DeclareMathOperator}
\nco{\halv}{half}
\nco{\order}{o} 
\nco{\ppower}{pp} 
\nco{\pcf}{pcf} 
\nco{\tcf}{tcf} 
\nco{\tlim}{tlim} 
\nco{\limtext}{lim} 
\nco{\prodt}{{\textstyle \prod}} 
\nco{\symdiff}{\triangle}
\nco{\dom}{dom}
\nco{\card}{card}
\nco{\lh}{lh}
\nco{\lgg}{lg}
\nco{\rge}{rge}
\nco{\otp}{otp}
\nco{\trunk}{tr}
\nco{\cf}{cf}
\nco{\nex}{next}
\nc{\uhr}{\restriction}
\nco{\supt}{supt}
\nco{\supp}{supp}
\nco{\Lim}{Lim}
\nco{\Leb}{Leb}
\nco{\modd}{mod}
\nco{\invariant}{inv}
\nco{\id}{id}
\nco{\RO}{RO}
\nco{\Dp}{Dp} 
\nco{\pss}{ps}
\nco{\acc}{acc}
\nco{\spec}{spec}
\nco{\pr}{pr}
\nco{\rt}{rt}
\nco{\suc}{suc}
\nco{\splitt}{split}

\nc{\potom}{\ensuremath{{\cal P}(\omega)}}
\nc{\potinf}{\ensuremath{[\omega]^\omega}}
\nc{\pfin}{\ensuremath{{\cal P}(\omega)/{\rm fin}}}

\nc{\potfin}{\ensuremath{[\omega]^{<\omega}}}
\nc{\inn}{\ensuremath{{\omega^{\uparrow \omega}}}}
\nc{\hoch}{^{<\omega}}
\nc{\hocho}{^{\omega}}
\nc{\tree}[1]{{[} #1 {]}_0}
\nc{\tre}[2]{ {#1}_{#2}}
   
\nc{\prooff}[1]{{\bf Proof} of #1:}
\nc{\proofend}{\makebox{} \hfill ${\bf \square}$ \\}
\nc{\proofendof}[1]{\makebox{} \hfill $\boldmath{\square}_{\rm #1}$ \\}
\nc{\beq}{\begin{eqnarray*}}
\nc{\eeq}{\end{eqnarray*}}
\nc{\bde}{\begin{list}}
\nc{\ede}{\end{list}}

\newenvironment{myrules}
{\begin{list}{}
{
 \setlength{\leftmargin}{0.5in}
 \setlength{\labelwidth}{1cm}
 \setlength{\labelsep}{0.2in}
 \setlength{\parsep}{0.5ex plus 0.2ex minus 0.1 ex}
 \setlength{\itemsep}{0.3ex plus 0.2 ex minus 0ex}
}}{\end{list}}

{\end{list}}

\newcounter{subalph}
{\end{list}}

\newcommand{\greek}[1]{\ifthenelse{\value{#1}=1}{\mbox{$\alpha$}}%
  {\ifthenelse{\value{#1}=2}{\mbox{$\beta$}}{%
   \ifthenelse{\value{#1}=3}{\mbox{$\gamma$}}{%
   \ifthenelse{\value{#1}=4}{\mbox{$\delta$}}{%
   \ifthenelse{\value{#1}=5}{\mbox{$\varepsilon$}}{%
   \ifthenelse{\value{#1}=6}{\mbox{$\zeta$}}{%
   \ifthenelse{\value{#1}=7}{\mbox{$\eta$}}{%
   \ifthenelse{\value{#1}=8}{\mbox{$\theta$}}{%
   \ifthenelse{\value{#1}=9}{\mbox{$\iota$}}{%
   \ifthenelse{\value{#1}=10}{\mbox{$\kappa$}}{%
   \ifthenelse{\value{#1}=11}{\mbox{$\lambda$}}{%
   \ifthenelse{\value{#1}=12}{\mbox{$\mu$}}{%
   \ifthenelse{\value{#1}=13}{\mbox{$\nu$}}{%
   \ifthenelse{\value{#1}=14}{\mbox{$\xi$}}{%
   \ifthenelse{\value{#1}=15}{\mbox{$\rm o$}}{%
   \ifthenelse{\value{#1}=16}{\mbox{$\pi$}}{%
   \ifthenelse{\value{#1}=17}{\mbox{$\varrho$}}{%
   \ifthenelse{\value{#1}=18}{\mbox{$\sigma$}}{%
   \ifthenelse{\value{#1}=19}{\mbox{$\tau$}}{%
   \ifthenelse{\value{#1}=20}{\mbox{$\upsilon$}}{%
   \ifthenelse{\value{#1}=21}{\mbox{$\varphi$}}{%
   \ifthenelse{\value{#1}=22}{\mbox{$\chi$}}{%
   \ifthenelse{\value{#1}=23}{\mbox{$\psi$}}{\mbox{$\omega$}%
  }}}}}}}}}}}}}}}}}}}}}}}}

\newcounter{subgreek}
{\end{list}}

\newcounter{subarabic}
{\end{list}}

\newcounter{subroman}
{\end{list}}

\newcount\skewfactor

\def\mathunderaccent#1#2 {\let\theaccent#1\skewfactor#2
\mathpalette\putaccentunder}
\def\putaccentunder#1#2{\oalign{$#1#2$\crcr\hidewidth
\vbox to.2ex{\hbox{$#1\skew\skewfactor\theaccent{}$}\vss}\hidewidth}}
\def\name{\mathunderaccent\tilde-3 }

\nc{\nname}{\name}

\nc{\even}{\ensuremath{\rm Even}}
\nc{\odd}{\ensuremath{\rm Odd}}

\nc{\al}{$\alpha$\  }
\nc{\om}{\omega}
\nc{\omm}{\ensuremath{\omega_1}}
\nc{\ep}{\varepsilon}
\nc{\tk}{\tilde{K}}
\nc{\concat}{\,\hat{} \,}   
\nc{\force}{\Vdash}
\nc{\fb}{f_{\bar{M}}}
\nc{\such}{\, : \,}   

\nc{\meager}{\ensuremath{{\cal M}}}
\nc{\lebesgue}{\ensuremath{{\cal N}}}
\nc{\nulll}{\ensuremath{{\cal N}}}
\nc{\ksigma}{\ensuremath{{\bf K}_\sigma}}
\nc{\ideal}{\ensuremath{{\cal I}}}
\nc{\ga}{\ensuremath{\frak a}}
\nc{\AAA}{{\cal A}}   
\nc{\gc}{\ensuremath{\frak c}}
\nc{\gs}{\ensuremath{\frak s}}
\nc{\gh}{\ensuremath{\frak h}}
\nc{\gd}{\ensuremath{\frak d}}
\nc{\gb}{\ensuremath{\frak b}}
\nc{\gro}{\ensuremath{\frak g}}
\nc{\gu}{\ensuremath{\frak u}} 
\nc{\gr}{\ensuremath{\frak r}} 
\nc{\gt}{\ensuremath{\frak t}}
\nc{\fff}{\ensuremath{\frak f}}
\nc{\gm}{\ensuremath{\mathfrak{mcf}}}
\nc{\gge}{\ensuremath{\mathfrak e}}
\nc{\cfupro}{\ensuremath{\cf(\upro)}}
\nc{\cfvpro}{\ensuremath{\cf(\vpro)}}
\nc{\gp}{\ensuremath{\frak p}}
\nc{\gk}{\ensuremath{\frak k}}

\nc{\add}[1]{\mbox{\ensuremath{{\rm add}(#1)}}}
\nc{\cov}[1]{\mbox{\ensuremath{{\rm cov}(#1)}}}
\nc{\unif}[1]{\mbox{\ensuremath{{\rm unif}(#1)}}}
\nc{\cof}[1]{{\mbox{\ensuremath{\rm cof}(#1)}}}

\nc{\addd}[2]{\mbox{\ensuremath{{\rm add}^{#1}(#2)}}}   
\nc{\covv}[2]{\mbox{\ensuremath{{\rm cov}^{#1}(#2)}}}   
\nc{\uniff}[2]{\mbox{\ensuremath{{\rm unif}^{#1}(#2)}}} 
\nc{\coff}[2]{{\mbox{\ensuremath{\rm cof}^{#1}(#2)}}}

\nc{\cd}{Cicho\'n's Diagram}
\nc{\COF}{\mbox{\bf Cof}}

\nc{\MA}{\mbox{\rm MA}}
\nc{\PFA}{\mbox{\rm PFA}}
\nc{\OCA}{\mbox{\rm OCA}}
\nc{\GCH}{\mbox{\rm GCH}}
\nc{\CH}{\mbox{\rm CH}}
\nc{\zfc}{\mbox{\rm ZFC}}
\nc{\sch}{\mbox{\rm SCH}} 
\nc{\ZF}{\mbox{\rm ZF}}
\nc{\NCF}{\mbox{\rm NCF}} 
\nc{\FD}{\mbox{\rm FD}}   

\nc{\fourG}{\mbox{\rm 4G}} 
\nc{\fourI}{\mbox{\rm 4I}}   

\nc{\Borelhood}{Borel measurability} 
\nc{\Pieinseins}{\mbox{${\bf \Pi}^1_1$}}
\nc{\seinseins}{\mbox{${\bf\Sigma}^1_1$}}
\nc{\seinszwei}{\mbox{${\bf\Sigma}^1_2$}}
\nc{\seinsdrei}{\mbox{${\bf\Sigma}^1_3$}}
\nc{\Deleinszwei}{\mbox{${\bf\Delta}^1_2$}}

\nc{\up}{\ensuremath{{\cal U}\mbox{\ensuremath{\rm -prod}}\,\omega}}
\nc{\upp}{\ensuremath{{\cal U}'\mbox{\ensuremath{\rm -prod}}\,\omega}}
\nc{\upro}{\ensuremath{{\cal U}\mbox{\ensuremath{\rm -prod}}\,\om}}
\nc{\fupro}{\ensuremath{f({\cal U})\mbox{\ensuremath{\rm -prod}}\,\om}}
\nc{\vpro}{\ensuremath{{\cal V}\mbox{\ensuremath{\rm -prod}}\,\om}}
\nc{\fpro}{\ensuremath{{\cal F}\mbox{\ensuremath{\rm -prod}}\,\om}}

\nc{\cff}[1]{{\text{cf}\,(#1)}}           
\nc{\cu}{\ensuremath{\cal U}}             
\nc{\ai}{\ensuremath{\forall^\infty}}     
\nc{\ei}{\ensuremath{\exists^\infty}}     
\nc{\ww}{\ensuremath{\omega^\omega}}      

\nc{\gw}{groupwise dense}

\nc{\kk}{car\-dinal cha\-rac\-teris\-tic}
\nc{\joker}{\ast}
\nc{\gtc}{Galois-Tukey connection} 

\nc{\av}[1]{{\rm Av}_{#1}}
\nc{\eps}{\varepsilon}
\nc{\n}{{\bf n}}                 
\nc{\m}{{\bf m}}

\nc{\marginparr}[1]{}
\nc{\footnoteee}{} 
\nc{\footnotee}{}  

\newcommand{\cal}{\mathcal}

\nc{\divs}{{c_0 \setminus \ell^1}}
\nc{\divser}{(\divs, \leq^*)/\thickapproy}
\nc{\bfin}{\RO(\pfin \setminus\{0\},\subseteq^*)}
\nc{\bdivser}{\RO(\divser)}
\nc{\inc}{{\rm INC}}
\nc{\com}{{\rm COM}}
\nc{\thickapproy}{\makebox{}\!\!\thickapprox}
\nc{\approy}{\makebox{}\!\!\approx}
\nc{\lessi}{\leqslant}
\nc{\gessi}{\geqslant}
\nc{\interior}[1]{{\rm int}(#1)}
\nc{\closure}[1]{{\rm cl}(#1)}
\nc{\Vo}{Vojt\'a\v{s}}

\nc{\precedeseq}{\leq^*} 
\nc{\precedes}{\prec}
\nc{\stronger}{\leqslant_{\bf P}}
\nc{\underlline}[1]{\hat{#1}}
\nc{\PO}{{\bf P}}
\nc{\charak}{\text{ch}}

\nc{\needed}{needed\ }
\nc{\neededc}{needed}
\nc{\Needed}{Needed\ }
\nc{\wneeded}{weakly needed\ }
\nc{\Wneeded}{Weakly needed\ }
\nc{\wneededc}{weakly needed}

\nc{\mup}{m_{\rm up}}
\nc{\mdn}{m_{\rm dn}}
\nco{\may}{may}
\nco{\aver}{av} 
\nco{\norm}{nor} 
\nco{\val}{val} 
\nco{\dis}{dis} 
\nco{\basis}{basis}
\nco{\pos}{pos}

\nc{\err}{\mbox{err}}
\nc{\eee}{\mbox{e}}
\nco{\Expect}{Exp}

\begin{document}


\title{Specializing Aronszajn trees by countable approximations}

\author{Heike Mildenberger and Saharon Shelah}


\thanks{The first author was partially supported by
a Minerva fellowship.}

\thanks{The second author's research 
was partially supported by the ``Israel Science
Foundation'', founded by the Israel Academy of Science and Humanities.
This is  the second author's work number 778.}

\address{Heike Mildenberger,
Institut f\"ur formale Logik, Universit\"at Wien, W\"ahringer Str.\ 25, A-1090 Wien, Austria}

\address{Saharon Shelah,
Institute of Mathematics,
The Hebrew University of Jerusalem,
 Givat Ram,
91904 Jerusalem, Israel
}

\email{heike@logic.univie.ac.at}
 
\email{shelah@math.huji.ac.il}

\begin{abstract}
We show that there are proper forcings based upon countable trees of creatures that 
specialize a given Aronszajn tree.
\end{abstract}



\maketitle


\setcounter{section}{-1}

\nothing{{\sf Comments and questions by Heike are in sans serif font.
Please take care of the questions written in this font.}}

\section{Introduction}\label{S0}

The main point of this work is finding
forcing notions specializing an Aronszajn tree, 
which are creature forcings, tree-like 
with halving, but being based on $\omega_1$ 
(the tree) rather than $\omega$.

\smallskip

For creature forcing in general there is ``the book
on creature forcing'' \cite{RoSh:470} and
for the uncountable case the work is extended in
\cite{RoSh:628} and \cite{RoSh:736}.
Since some of the main premises made in the mentioned work
are not fulfilled in our
setting, it  serves mainly as a guideline, whereas
numerous technical details here are different and new.

\nothing{shall give an almost self-contained presentation,
referring to citations only if they fit almost exactly.
}
\smallskip

The norm of creatures (see Definition~\ref{1.7}) we shall use 
is natural for specializing Aronszajn
trees. It is convenient if there is some $\alpha < \omega_1$ such that
the union of the domains of the partial specialization functions
that are attached to any branch of the tree-like forcing condition
 is the initial segment of the Aronszajn tree $(T_A)_{<\alpha}$,
i.e.\ the union of the levels less than $\alpha$. However,  allowing for every branch
of a given condition finitely many
possibilities $(T_A)_{<\alpha_i}$ with finite sets $u_i$ sticking
out of $(T_A)_{<\alpha_i}$
 is used for
density arguments that show that the generic filter
leads to a total specialization function.

\nothing{{\sf Saharon, you added:}
Relation to [RoSh:470]?
{\sf Saharon, what do you want to say here?}
}
\section{Tree creatures}\label{S1}

In this section we define the tree creatures  which
will be used later to describe the branching of the countable
trees that will serve as forcing conditions.
We prove three important technical properties about gluing together (Claim~\ref{1.9}),
about filling up (Claim~\ref{1.10}) and about changing the base together with thinning out
(Claim~\ref{1.11}) of creatures.
We shall define the forcing conditions only in the next section. They will be countable
trees with finite branching, such that each node and its immediate successors
in the tree are described by a creature in the sense of Definition~\ref{1.5}.
Roughly spoken, in our context, a creature will be an arrangement of partial 
specialization functions with some side conditions. 

\smallskip

We reserve the symbol $(T,\triangleleft_T)$ for the trees in the forcing
conditions, which are trees of partial specializiation functions 
of some given Aronszajn tree $(T_A,<_{T_A})$. 
A specialization function is a function
$f \colon T_A \to \omega$ such that
for all $s,t \in T_A$, if $s <_{T_A} t$, then $f(s) \neq f(t)$,
see \cite[p.\ 244]{Jech}.

$\chi$ stands for some sufficiently high regular cardinal, and ${\mathcal H}(\chi)$ denotes
the set of all sets of hereditary cardinality less than $\chi$. For our purpose, $\chi =
(2^\omega)^+$ is enough.

\smallskip

Throughout this work we make the following assumption:

\begin{hypothesis}\label{1.1}
$T_A$ is an Aronszajn tree ordered by $<_{T_A}$, and 
for $\alpha < \omega_1$ the level $\alpha$ of $T_A$ satisfies:
$$(T_A)_\alpha \subseteq [\omega \alpha, \omega \alpha + \omega).$$
\end{hypothesis}

Throughout this work, $T_A$ will be fixed.
We define the following finite approximations of
specialization maps:

\begin{definition}\label{1.2}
For $u \subseteq T_A$ and $n< \omega$ we let
$$\spec_n(u) = 
\{ \eta \mid \eta \colon u \to [0,n) \; \wedge \; (\eta(x) = \eta(y)
\rightarrow \neg(x <_{T_A} y)) \}.$$
We let $\spec(u)= \bigcup_{n<\omega} \spec_n(u)$ and
$\spec = \spec^{T_A} =
\bigcup\{ \spec(u) \such u \subset T_A, u \mbox{ finite}\}.$  
\end{definition}

\begin{choice}\label{1.3}
We choose three sequences of natural numbers $
\langle n_{k,i} \such i < \omega \rangle$,
$k=1,2,3$,
such that the following growth conditions are fulfilled:
\begin{eqnarray}
i \cdot n_{1,i} &<& n_{3,i},\\
n_{2,i} &<& n_{1,i+1},\\
n_{1,i} \cdot n_{1,i} &\leq& n_{1,i+1}, \\
n_{1,i} &\leq & n_{2,i}.
\end{eqnarray}

We fix them
for the rest of this work.
\end{choice}

\smallskip

\nothing{
{\sf Maybe we should switch $n_{3,i}$ and $n_{2,i}$ because of the second
equation. The numbers $n_{0,i}$ are so far not used at all in this work.}
}

\smallskip

We compare with the book \cite{RoSh:470} in order to justify
the use of he name creature. However,  we cannot just cite that work, because the
framework developed there is not suitable for the approximation of uncountable domains $T_A$.

\begin{definition}\label{1.4} 
\begin{myrules}
\item[(1.)]
\cite[1.1.1]{RoSh:470} 
A triple $t=(\norm[t],\val[t],\dis[t])$ is a 
a weak creature for ${\bf H}$ if 
\begin{myrules}
\item[(a)] $\norm[t] \in {\mathbb R}^{\geq 0}$,
\item[(b)] Let ${\bf H} = \bigcup_{i \in \omega} {\bf H}(i)$
and let ${\bf H}(i)$ be sets. Let $\triangleleft$ be the strict
initial segment relation.

$\val[t]$ is a non-empty subset of 
$$\left\{ \langle x,y \rangle \in \bigcup_{m_0 < m_1 < \omega} [
\prod_{i< m_0}{\bf H}(i) \times \prod_{i<m_1} {\bf H}(i)]
\such x \triangleleft y \right\}.$$

\item[(c)] $\dis[t] \in {\mathcal H}(\chi)$.
\end{myrules}

\item[(2.)]
$\norm$ stands for norm, $\val$ stands for value, and $\dis$ stands
for distinguish.
\end{myrules}
\end{definition}

In our case, we drop the component $\dis$ (in the case of simple creatures in the sense of
Definition~\ref{1.5})
or it will be called $k$ (in the case of creatures), 
an additional coordinate, which is a natural 
number.
In order to stress some parts of the weak creatures
$t$ more than others,
we shall write $\val[t]$ in a slightly different form and call it 
a simple creature, $\bf c$.

\smallskip

As we will see in the next definition, in this work
(b) of \ref{1.4} is not fulfilled: For us
$\val$ is a non-empty subset of
$\{ \langle x,y\rangle \in \spec \times \spec \such x \triangleleft_T y \}$
for some partial order $\triangleleft_T$ as in Definition~\ref{2.1}.
Though the members of $\spec$ are finite partial functions, they
cannot be written with some $n \in \omega$ as a domain, since
$\spec$ is uncountable and we want to allow arbitary finite parts.
Often properness of a tree creature forcing follows
from the countability of ${\bf H}$. Note that our analogue
to ${\bf H}$ is
not countable. In Section~3 we shall prove that the notions of forcing
we introduce are proper for other reasons.

\smallskip

Nevertheless the simple creature in
the next definition is a specific case for
the value of a weak creature in the sense of \ref{1.4} without item (1.)(b), 
and the creature from
the next definition can be seen as a case of a value and
a distinction part of a weak creature.

\begin{definition}\label{1.5}
\begin{myrules}

\item[(1)]
A simple creature 
is a tuple ${\bf c} =(i({\bf c}),\eta({\bf c}), \rge(\val({\bf c})))$
with the following properties:

\begin{myrules}
\item[(a)]
The first component, $i({\bf c})$, is called the kind of ${\bf c}$ and is just a natural number.

\item[(b)] The second component, $\eta({\bf c})$, is called the base of $\bf c$.
We require ($\eta({\bf c}) = \emptyset$ 
and
$i({\bf c}) = 0$) or ($i({\bf c})$ is the smallest $i$ such that
$|\dom(\eta({\bf c}))| \leq n_{2,i-1}$), and $\eta({\bf c}) \in \spec_{n_{3,i-1}}$.

\item[(c)] The range of the value of $\bf c$,
$\rge(\val({\bf c}))$, is a non-empty subset of
$\{ \eta \in \spec_{n_{3,i}} \such \eta({\bf c}) 
\subseteq \eta \wedge |\dom(\eta)| < n_{2,i}\}$,
such that
$|\rge(\val({\bf c}))| < n_{1,i}$.

\nothing{{\sf There is ellbow room: Finite instead of  bounded by 
$n_{1,i}$ would be enough so far, for all the applications of
K\"onig's lemma.}}

\nothing{{\sf In the manuscript, it is always $n_{3,i}$. I changed them to
$n_{2,i}$ and $n_{1,i}$. Maybe it should be vice versa. I did not
use some order relation between $n_{k,i}$ for different $k$'s in the proofs.
I wrote ``$\rge(\val)$'' instead of ``$\val$'' 
in accordance with \cite{RoSh:470}.}}

So we have $\val({\bf c})= \{\eta({\bf c})\} \times \rge(\val({\bf c}))$.
That the domain is a singleton, is typical for tree-creating
creatures.

\item[(d)] If $\eta_1 \in \rge(\val({\bf c}))$ 
and  $x
\in \dom(\eta_1) \setminus \dom(\eta({\bf c}))$
then there is some $\eta_2 \in \rge(\val({\bf c})$ such that
$x \in \dom(\eta_2) \rightarrow
\eta_1(x) \neq \eta_2(x)$.
\end{myrules}

\item[(2)] A creature ${\bf c}^+ $ is 
a tuple $(i({\bf c}^+),\eta({\bf c}^+), \rge(\val({\bf c}^+)),k({\bf c}^+))$
where \\
$(i({\bf c}^+),\eta({\bf c}^+), \rge(\val({\bf c}^+)))$
is a simple creature, and $k({\bf c}^+)\in \omega$ is an additional
coordinate.

\item[(3)] An  (simple) $i$-creature is a (simple) creature
with $i({\bf c}^+ ) = i$ ($i({\bf c}) = i$).

\item[(4)] If ${\bf c}^+$ is a creature we mean by $\bf c$ the
simple creature such that ${\bf c}^+ =({\bf c}, k({\bf c}^+))$.

\item[(5)] The set of creatures is denoted by $K^+$, 
and the set of simple creatures is denoted by $K$.

\end{myrules}

\end{definition}

\begin{remark}\label{1.6}
 By \ref{1.5}(d) we have that $\eta({\bf c}) = \bigcap\{
\eta \such \eta \in \rge(\val({\bf c}))\}$,
\nothing{{\sf Saharon does not think that this follows, but Heike thinks
that this follows from (d)}}
and also $i({\bf c})$ 
is determined by
$\eta({\bf c})$ and hence from $\rge(\val({\bf c}))$. 
Thus, in our specific case, every simple creature is
determined by the range of its value.
\end{remark}

For a real number $r$ we let $m=\lceil r \rceil$ be the smallest 
natural number such that $m \geq r$. So, for negative numbers
$r$, $\lceil r \rceil =0$. We let
$\lg$ denote the logarithm function to the base 2.
Let $\log_2(x) = \lceil \lg(x) \rceil$ for $x > 0$, and
we set $\log_2 0 =0$.

\nothing{For $f \colon \omega \to \omega$ we set $f^{(0)} = f$ and 
$f^{(k+1)} = f\circ f^{(k)}$ for $k \in \omega$.
For a natural number $n >0$ we let $\log_\ast n $ be the smallest
 $k$ such that
$\log_2^{(k)} n =0$.
{\sf Warning. These functions are not used in the paper}
}

\begin{definition}\label{1.7}
\begin{myrules}
\item[(1)] For a 
simple $i$-creature ${\bf c}$ 
we define $\norm^0({\bf c})$ as the maximal natural number $k$ such that
if $a \subseteq n_{3,i}$ and $|a| \leq k$  and $B_0, \dots , B_{k-1}$ are 
branches of $T$, then there is $\eta \in \val({\bf c})$ such that
\begin{myrules}
\item[($\alpha$)] $(\forall x \in 
(\bigcup_{\ell < k} B_\ell \cap \dom(\eta))) \setminus
\dom(\eta({\bf c})))(\eta(x) \not\in a)$,
\item[($\beta$)] $\displaystyle{\frac{|\dom(\eta)|}{n_{2,i}} 
\leq \frac{1}{2^k}}$. 
\nothing{{\sf I changed the $n_{3,i}$ in the denominator
to $n_{2,i}$, according to \ref{1.4}(c).}
\nothing{I think we should add $\displaystyle{\frac{|\rge(\eta)|}{n_{3,i}} 
\leq \frac{1}{2^k}}$ as well}
}
\end{myrules}

\item[(2)]
\nothing{We consider two versions for
$\norm^{\frac{1}{2}}({\bf c})$.\\
{\em Version 1:} $\norm^{\frac{1}{2}}({\bf c})=\norm^0({\bf c})$.\\
{\em Version 2:} }
We let $\norm^\ast({\bf c}) =
\log_2(\frac{n_{1,i({\bf c})}}{|\val({\bf c})|}))$, and
\nothing{{\sf I changed $n_{3,i}$ to $n_{1,i}$ here.}}
 $\norm^{\frac{1}{2}}({\bf c})=\min(\norm^0({\bf c}),\norm^\ast({\bf c}))$.
\nothing{Version 2 has some stronger properties.
}

\item[(3)] 
We define
$\norm^1({\bf c}) = 
\log_2(\norm^{0}({\bf c}))$,
and $\norm^2({\bf c}) = \log_2(\norm^{\frac{1}{2}}({\bf c}))$.

\nothing{{\sf So far, only the first version of
$\norm^1$ is used in the paper. What do you intend to prove
with this variation? I strongly recommend that we give
it a different name.}}

\nothing{When using $\norm^1$ we shall always indicate 
which version was used.
{\sf I changed the $n_{3,i}$ in the denominator to
$n_{0,i}$. This is the only appearance of $n_{0,i}$.}
}
\item[(4)] 
In order not to fall into specific computations, we use
functions $f$ that exhibit the following properties, in order
to define norms on (non-simple) creatures:

\begin{myrules}
\item[$(\ast)_1$] $f$ is a two-place function.
\item[$(\ast)_2$] $f$ fulfils the following monotonicity properties:
If $n_1 \geq n_2 \geq k_2 \geq k_1$ then
$f(n_1,k_1) \geq f(n_2,k_2)$.
\item[$(\ast)_3$] (For the 2-bigness, see Definition~\ref{1.12})
$f(\frac{n}{2},k) \geq f(n,k) -1$.
\item[$(\ast)_4$] $n \leq k \rightarrow f(n,k)=0$.
\item[$(\ast)_5$] (For the halving property, see Definition~\ref{3.3})
For all $n,k,i$: If $f(n,k) \geq i+1$, then there is some $k'(n,k)=k'$ 
such that
$k<k'<n$ and for all $n'$, if $k'<n'<n$, then
$$f(n',k) = f(n',k') + f(k',k) \geq f(n,k)-1.$$
\nothing{{\sf I changed point 5 considerably in comparison to what we had on March 7 on the
yellow pages, in order to be able to prove \ref{3.5}}}
\end{myrules}
\end{myrules}
\end{definition}

For example, $f(n,k) = \lg(\frac{n}{k})$ and
$k'(n,k) = \sqrt{nk}$, fulfil these conditions.
For a  creature ${\bf c}^+$ we define
its norm
 $$\norm({\bf c}^+) = f(\norm^{\frac{1}{2}}({\bf c}),k({\bf c}^+)).$$

\nothing{
where $f\colon {\mathbb R} \to {\mathbb R}$ is monotonic, and
fulfils for all $x,y$, $f(x\cdot y) \geq f(x) + f(y)$ and
$f(\sqrt{x}) \geq f(x) -1$.}

\begin{remark}\label{1.8}
1. Note that property (1)(d) of simple creatures (Definition~\ref{1.5}) follows from
$\norm^0({\bf c}) >0$. So we will not check this property any more, but restrict
ourselves to creatures with strictly positive $\norm^0$. 

2. Definition~\ref{1.7} speaks about infinitely many requirements, by ranging over all
$k$-tuples of branches of $T$. However, at a crucial point in the
proof of Claim~\ref{1.10}
 this boils down to
counting the possibilities for $a$.
\end{remark}
 
\nothing{
{\sf Saharon, I changed this 
as in the new yellow pages from November 2000.
Originally (in the yellow pages and also in Mirna's Typescript) there
was something else with $\log_\ast$.}}
\nothing{I did not use it at all, because I did not
see how to get halving in the sense of 
\cite[Definition 2.2.7]{RoSh:470}. I also do not see at all
how halving could help in \ref{3.9} of this work}

\smallskip

The next claim shows that we can to extend the functions in the value of a
creature and at the same
time decrease the norm of the creature only  by a small amount.

\begin{claim}\label{1.9}
Assume that 
\begin{myrules}
\item[(a)]$\eta^\ast \in \spec$,

\item[(b)] $\bf c$ is a simple $i$-creature with base $\eta^\ast$, 
$\norm^{0}({\bf c}) >0$,

\item[(c)] $k^\ast > 1$, $|\rge(\val({\bf c}))| \cdot k^\ast \leq n_{1,i}$,

\item[(d)] for each $\eta \in \rge(\val({\bf c}))$ and $k <k^\ast$ we are given
$\eta \subseteq \rho_{\eta,k} \in \spec_{n_{3,i}}$
with $|\dom(\rho_{\eta,k})| < n_{2,i}$, 

\item[(e)] for each $\eta \in \rge(\val({\bf c}))$, 
if $k_1 < k_2 < k^\ast$ and $x_1 \in \dom(\rho_{\eta,k_1})
\setminus \dom(\eta)$ and
$x_2 \in \dom(\rho_{\eta,k_2}) \setminus \dom(\eta)$, then 
$x_1$, $x_2$ are  $<_{T_A}$-incomparable,
\nothing{{\bf  With this strong premise, you can never fill up
the union of the domains of the members of spec along a branch in
a condition to some $T_{<\alpha}$. What is the claim good for?}}

\item[(f)] $\ell^\ast=
\max\{|\dom(\rho_{\eta,k})| +1 \such \eta \in \rge(\val({\bf c})) 
\wedge k < k^\ast\}$,
\nothing{\item[(f)] For $\eta \in \val({\bf c})$, $k<k^\ast$ we find
$\langle m_{\eta,k,j}  \such j < k^\ast
\rangle$ with no repetition such that
\begin{myrules}
\item[($\alpha$)] $m_{\eta,k,j} < n_{3,i}$ for all $k,j < k^\ast$,
\item[($\beta$)]
$\rho_{\eta,k} \cup
\{ (y_\eta,m_{\eta,k,j})\} \in \spec$ for all $k,j < k^\ast$,0
\end{myrules}}

\end{myrules}

Then 
there is a simple $i$-creature $\bf d$ given by
$$\rge(\val({\bf d}))=
\{ \rho_{\eta, k}  \such k < k^*,
\eta \in \rge(\val({\bf c}))\}.$$
We have $\eta({\bf d}) = \eta^\ast$, and  $\norm^0({\bf d}) \geq m_0 \stackrel{\rm{\small def}}{=}
\min\left\{\norm^0({\bf c}), \log_2(\frac{n_{2,i}}{\ell^*}), 
k^\ast -1\right\}$.

\end{claim}

\proof First of all we are to check Definition~\ref{1.5}(1).
Clauses (a),(b), and (c) follow immediately from the premises of the claim.
${\bf d}$ satisfies clause (d): This follows from the proof of the inequality for
$\norm^0({\bf d})$ below.
\relax From premise (e) and from the properties of ${\bf c}$ it follows that $\eta({\bf d}) = \eta^\ast$.

\smallskip

Now for the norm:
We check clause $(\alpha)$ of Definition~\ref{1.7}.
Let branches $B_0, \dots B_{m_0-1}$ of $T_A$ and
a set $a \subseteq n_{3,i}$ be given, $|a| \leq m_0$.
Since $m_0 \leq \norm^0({\bf c})$, there is some $\eta \in 
\rge(\val({\bf c}))$ such that
$(\forall x \in \left(\bigcup_{\ell <  m_0} B_\ell\right) \cap \dom(\eta) 
\setminus \dom(\eta({\bf c})))(\eta(x) \not\in a)$.
We fix such an $\eta$. Now for each $\ell < m_0$, we let
$$w_{\eta,\ell}= \{ j < k^\ast \such \exists x \in B_\ell \cap
\dom(\rho_{\eta,j}) \setminus \dom(\eta) \}.$$

Now we have that $|w_{\eta,\ell}| \leq 1$ because otherwise
we would have $k_1 < k_2 < k^\ast$ in $w_{\eta,\ell}$ and
$x_i \in B_\ell \cap \dom(\rho_{\eta,k_i}) \setminus \dom(\eta)$.
As $x_1 $ and $x_2$ are $<_{T_A}$-comparable, this is contradicting 
the requirement $(d)$ of \ref{1.9}.

Since $m_0 < k^\ast$, there is some $j \in k^\ast \setminus
\bigcup_{\ell < m_0} w_{\eta,\ell}$. For such a $j$,
$\rho_{\eta,j}$ is as required.

\smallskip

We check clause $(\beta)$ of Definition~\ref{1.7}.
We take the $\rho_{\eta,j}$ as chosen above. Then we have 
$$\frac{|\dom(\rho_{\eta,j})|}{n_{2,i}}\leq \frac{\ell^*}{n_{2,i}}
= \frac{1}{2^{\log_2\left(\frac{n_{2,i}}{\ell^*}\right)}} \leq
\frac{1}{2^{m_0}},$$
as $m_0 \leq \log_2\left(\frac{n_{2,i}}{\ell^*}\right)$.
\proofend

Whereas the previous claim will be used only in Section 3
in the proof on properness (see Claim~\ref{3.8}),
the following two claims will be used in the next section for
density arguments in the forcings built from creatures.

\begin{claim}\label{1.10}
\nothing{{\sf 10 A on page D in amendment,
but I changed it for $m$ $x$'s at the same time,  because this is 
used later, and it not just $m$ times the application of
the claim because of premise (d).}}
Assume
\begin{myrules}

\item[(a)] $\bf c$ is a simple $i$-creature.
\item[(b)] $k=\norm^0({\bf c}) \geq 1$ and $k \leq n_{1,i}$.
\item[(c)] $x_0, \dots x_{m-1} \in T$, $1\leq m \leq 
\min(k, \frac{n_{2,i}}{2^k})$.
\item[(d)]  $|\rge(\val({\bf c}))| \cdot \binom{k}{m} \leq n_{1,i}$.
\item[(e)] If $\eta \in \rge(\val({\bf c}))$, then
$|\{y \in \dom(\eta) \such (\exists m' < m) (x_{m'} <_T y) \}|<i$.

\end{myrules}

Then there is ${\bf d}$ such that $\eta({\bf d})
=\eta({\bf c})$ and
\begin{multline*}
\rge(\val({\bf d})) \subseteq \{ \nu \in \spec_{3,i} \such |\dom(\nu)| \leq n_{2,i},\\
(\exists\eta \in \rge(\val({\bf c})))( \eta \subseteq \nu \wedge
 \dom(\nu) = \dom(\eta) \cup\{x_0,\dots x_{m-1}\}) \},
\end{multline*}
 
such that $|\rge(\val({\bf d}))| \leq n_{1,i}$ and such that
for each $\eta$ there are sufficiently many ($\binom{k}{m}$ suitable ones
 suffice) extensions $\nu \in \rge(\val({\bf d}))$
as described in the first paragraph of the proof.

Moreover we can choose $\bf d$ such that 
\begin{myrules}
\item[$(\alpha)$] $\bf d$ is a simple $i$-creature.
\item[$(\beta)$] $\norm^0({\bf d}) \geq k-m$.
\nothing{\item[$(\gamma)$] For $\nu 
\in \val({\bf d})$ we have that $x \in \dom(\nu)$.}
\end{myrules}
\end{claim}

\proof
\nothing{Remark: We may add some $x$'s, this does not matter, but in fact
we add them one by one.} 
\nothing{
In order to see that $\rge(\val({\bf d}))
\neq \emptyset$, we first
take some $\eta \in \rge(\val({\bf c}))$ such that
$|\dom(\eta)| \leq \frac{n_{2,i}}{2^k}$.
Then we take $\nu=\eta \cup \{(x,z)\}$.
take $\nu(x) \neq \eta(y)$ for
$y \in \dom(\eta)$, $x$ $<_{T_A}$-comparable to $y$. Since $i<n_{3,i}$
and since for $\eta \in \rge(\val({\bf c})) $ $|\dom(\eta)| \leq
\frac{n_{1,i}}{2^k}$, this is possible.
}

We take for $m' < m$,  $z_{m'} \in n_{3,i}
\setminus (\bigcup_{\eta \in \rge(\val({\bf c}))} \{ 
\eta(y) \such x <_{T_A} y \} \cup \{ z_{m''} \such m'' < m'\})$. 
Since $|\rge(\val({\bf c}))|<n_{1,i}$ and
by (d) and since by \eqref{1.1}
$(i-1) \cdot n_{1,i}  +k < n_{3,i}$ there is such a $z_{m'}$,
and indeed, which is important for ${\bf d}$ being a creature and for its norm, there are at least
$k-m'$ such $z_{m'}$'s for every $\eta\in \rge(\val({\bf c}))$.
We take  all these choices
$\nu_{\eta,\bar{z}} = \eta \cup \{(x_{m'},z_{m'}) \such m' < m \}$ into $\rge(\val({\bf d}))$.
Hence we  can choose all $\nu_{\eta,\bar{z}}$
so that we avoid any given $a$ of size $k-m$ with all the $z_{m'}$'s.
\smallskip

Now we check the norm:
Let $B_0,\dots,B_{k-m-1}$ be branches of $T_A$ and let
$a \subseteq n_{3,i({\bf c})}$, $|a| \leq k-m$. 
We have to find $\nu \in \rge(\val({\bf d }))$ such that 
$(\forall \ell < k-m)(\forall y \in \dom(\nu) 
\cap B_\ell \setminus \dom(\eta({\bf c})) (\nu(y) \not\in a)$ and
$|\dom(\nu)| \leq \frac{n_{2,i}}{2^{k-m}}$.
For $m'<m$ we choose
$B_{k-m-1+m'}$, a branch containing $x_{m'}$. 
We take for $m' < m$, $z_{m'} \in n_{3,i}
\setminus (\bigcup_{\eta \in \rge(\val({\bf c}))} \{ 
\eta(y) \such x <_{T_A} y \} \cup a \cup\{z_{m''}\such m'' < m'\})$.
We set $a' = a \cup \{z_0,\dots , z_{m-1}\}$.

By premise (b), we find
 $\eta \in \rge(\val({\bf c}))$ for $a'$
and $B_0, \dots , B_{k-1}$ such that
\begin{eqnarray*}
(1)&& (\forall \ell < k-1)(\forall x \in \dom(\eta) 
\cap B_k \setminus \dom(\eta({\bf c})) (\eta(x) \not\in a')
\mbox{ and} \\
(2)&&|\dom(\eta)|  \leq \frac{n_{2,i}}{2^k}.
\end{eqnarray*}

Now $\nu_{\eta,\bar{z}}$  is a witness for the norm.
We have $\frac{n_{2,i}}{2^k} +m \leq \frac{n_{2,i}}{2^{k-m}}$, which 
follows from the premises on $m$.
The only thing to show is that $\nu$ is really a specialization function.
So let $y \in \dom(\eta)$ and $y<_{T_A} x_{m'}$. Then
$\nu(y)=\eta(y) \neq \nu(x_{m'})= z_{m'}$, because 
$y$ is on the branch leading to
$x_{m'}$ and because of (1). If $y>_{T_A} x_{m'}$, then we have 
taken care of $y$ simultaneously for 
all $\eta$'s by our choice of the $z_{m'}$'s.
\proofend

An analogous version of Claim~\ref{1.10} with $\norm^{\frac{1}{2}}$ instead of
$\norm^0$ holds as well. The analogous requirements to premises (c) and (d) are
even easier: If we work with $\norm^{\frac{1}{2}}$ and use
$n_{1,i} \leq n_{2,i}$ from equation (1.4) in the Choice~\ref{1.3}, 
then $1 \leq m \leq k$ is enough in premise  (c).
Premise (d) is included in $\norm^{\frac{1}{2}}({\bf c}) = k$.

\smallskip

Suppose we have filled up the range of the value
of a creature according to
one of the previous claims. Then we want that these
extended functions can serve as bases for suitable creatures as
well. This is provided by the next claim.

\begin{claim}\label{1.11}
\nothing{(2.7A on page D$_\alpha$ and D$_\beta$)} 
Assume that
\begin{myrules}
\item[(a)]
${\bf c}$ is a simple $i$-creature.
\item[(b)]
$k=\norm^0({\bf c}) \geq 1$.
\item[(c)]
$\eta^\ast \supseteq \eta({\bf c})$,
$\eta^\ast \in \spec_{n_{3,i}}$ (note that
we do not suppose that $\eta^\ast \in
\rge(\val({\bf c}))$). 
Furthermore we assume $|\dom(\eta^\ast)| \leq n_{2,i({\bf c})-1}$.
\item[(d)]
We set
$$\ell^\ast_2=|\dom(\eta^\ast) 
\setminus \dom(\eta({\bf c}))|,$$ and 
\begin{equation*}
\begin{split}
\ell_1^\ast= |\{y \such & (\exists \nu
\in \rge(\val({\bf c})))(y \in \dom(\nu) \setminus \dom(\eta({\bf c})))\\
& \wedge (\exists x \in \dom(\eta^\ast) 
\setminus \dom(\eta({\bf c})))(x<_{T_A}y) \}|.
\end{split}
\end{equation*}
and we assume that $\ell^\ast_1 + \ell_2^\ast < \norm^0({\bf c})$.
\end{myrules}
We define
${\bf d}$ by  $\eta({\bf d}) = 
\eta^\ast$ and $$\rge(\val({\bf d})) =
\{ \nu \cup \eta^\ast \such \nu \in \rge(\val({\bf c}))
\; \wedge \; 
\nu \cup \eta^\ast \in \spec_{n_{3,i}}, |\dom(\nu \cup \eta^\ast)| \leq n_{2,i}\}.$$
Then \begin{myrules}
\item[$(\alpha)$] ${\bf d}$ is a simple $i$-creature.
\item[$(\beta)$] 
$\norm^0({\bf d}) \geq \norm^0({\bf c}) - 
\ell_2^\ast - \ell_1^\ast$.
\end{myrules}
\end{claim}

\proof Item $(\alpha)$ follows from the requirements on
$\eta^\ast$ and from the estimates on the norm, see below. For item $(\beta)$, we set
$k = \norm^0({\bf c}) - \ell_1^\ast - \ell_2^\ast$.
 We let $B_0, \dots B_{k-1}$ be branches of $T_A$ and $a \subseteq 
n_{3,i({\bf c})}$, $|a| \leq k$. We set $\ell^\ast =
\ell_1^\ast + \ell_2^\ast$. We let $\langle y_\ell 
\such 
\ell < \ell_1^\ast \rangle$ list 
$Y=\{ y \such \exists \nu ( \nu \in \rge(\val({\bf c})) \wedge
y \in \dom(\nu)) \wedge
\exists x(x\in \dom(\eta^\ast) \setminus \dom(\eta({\bf c})) \wedge x\leq_{T_A} y )) \}$ without repetition.
Let $B_k,\dots,B_{k+\ell^\ast -1}$ be branches of $T_A$ such that $y_\ell 
\in B_{k+\ell}$ for $\ell < \ell^\ast_1$. Let $\langle
x_\ell \such \ell < \ell_2^\ast\rangle$ list $\dom(\eta^\ast) \setminus 
\dom(\eta({\bf c}))$. Take for $\ell < \ell^\ast_2$,
 $B_{k+\ell_1^\ast+\ell}$ such that
$x_\ell \in B_{k + \ell_1^\ast + \ell}$. We set $a'=a \cup
\{\eta^\ast(x_\ell) \such \ell < \ell^\ast_2 \}$.
Since $\norm^0({\bf c}) \geq k + \ell^\ast$ there is 
some $\nu \in \rge(\val({\bf c}))$ such 
that $\forall x \in ((\dom(\nu) \setminus 
\dom(\eta({\bf c}))) \cap \bigcup_{\ell < k + \ell^\ast} B_\ell ) 
(\nu(x) \not\in a')$.
Then, if $x \not\in \dom(\eta^\ast)$, $(\nu \cup \eta^\ast)(x) \not\in a$. 
Moreover $|\dom(\nu\cup \eta^\ast)| \leq \frac{n_{2,i}}{2^k} + \ell_2^\ast
\leq \frac{n_{2,i}}{2^{k-\ell}}$, if $\frac{n_{2,i}}{2^k}$ is large enough.
(This premise will always be fulfilled in our applications, because
$n_{1,i} \leq n_{2,i}$. We just
perform all our operations on forcing conditions only at high levels $i$, 
compared to the
size of the given $\eta^\ast$.
This will be done in the next section.)

We have to show that 
that if $\nu \cup \eta^\ast$ is a partial specialization:
Since $\eta^\ast$ and $\nu$ are specialization maps,
we have to consider only the case $x \in \dom(\eta^\ast)
\setminus \dom(\eta({\bf c}))$ and ($y \in Y$ or ($y \in 
\dom(\nu) \setminus \dom(\eta^\ast))$ and $y<_{T_A}x$)).
If $y \in Y$, then we have $\nu(y) \neq \eta^\ast(x_\ell)$ for all
$\ell < \ell^\ast_2$. If $y \in \dom(\nu) \setminus \dom(\eta^\ast)$ and
$y<_{T_A} x$, then $y$ is in a branch leading to some
$x_\ell$, $\ell< \ell^\ast_2$, and hence again $\nu(y) \neq
\eta^\ast(x_\ell)$, $\ell<\ell^\ast_2$.
\proofend

\bigskip

In the applications, the proofs of 
the density properties, $\ell^\ast_2$ will
be small compared to the norm (we add $\ell_2^\ast$ points to the domain
of the functions in the range of the value of a creature with
sufficiently high norm)
and $\ell_1^\ast \leq |u|$, were $u$ is the set that sticks out of $(T_A)_{<
\alpha(p)}$ (see Definition~\ref{2.2} and Remark~\ref{2.5}).
 We will suppose that these two are small
in comparison to $\norm^0({\bf c})$, so that the premises for Claim~\ref{1.11} are
fulfilled.

\smallskip

We also need Claims 1.9, 1.10 and 1.11 for $\norm^{\frac{1}{2}}$ instead of $\norm^0$.
This is proved by easy but a bit tedious accounting of $\norm^\ast({\bf c}) = 
\log_2(\frac{n_{1,i({\bf c})}}{\val({\bf c})})$.
Just see that $|\val({\bf c})|$ increases 
only in a controllable way in  Claim~\ref{1.9} and in Claim~\ref{1.10}
and does not increase at all in Claim~\ref{1.11}.
Hence also if $\norm^\ast$ is the part determining the minimum in 
$\norm^{\frac{1}{2}}$,
the latter falls at most by $k^\ast$ in \ref{1.9} from ${\bf c}$ to ${\bf d}$,
and at most by $\log_2(\binom{k}{m})$ in \ref{1.10} and does not decrease in \ref{1.1}.

\bigskip

The next claim will help to find large homogeneous subtrees of 
the trees built from creatures that will later be used a forcing conditions.
\nothing{As we shall see in the next section,
this property has useful consequences for the tree-like forcing
notions built from the tree-creatures.
}

\begin{claim}\label{1.12}
\begin{myrules}
\item[(1)] The 2-bigness property 
\cite[Definition 2.3.2]{RoSh:470}. If ${\bf c}$ is a simple
$i$-creature with $\norm^1({\bf c}) \geq k+1$, and 
${\bf c}_1$, ${\bf c}_2$ are simple $i$-creatures
such that $\val({\bf c})=\val({\bf c}_1) \cup
\val({\bf c}_2)$, then $\norm^1({\bf c}_1) \geq k$ 
or $\norm^1({\bf c}_2) \geq k$.
The same holds for $\norm^2$.

\item[(2)] 
 If ${\bf c}^+$ is a 
$i$-creature with $\norm({\bf c}) \geq k+1$, and 
${\bf c}^+_1$, ${\bf c}^+_2$ are  $i$-creatures
such that $\val({\bf c})=\val({\bf c}_1) \cup
\val({\bf c}_2)$, and $k({\bf c}^+_1)=k({\bf c}^+_2)=k({\bf c}^+)$, 
then $\norm({\bf c}_1^+) \geq k$ 
or $\norm({\bf c}^+_2) \geq k$.
\end{myrules}
\end{claim}

\proof
(1) We first consider $\norm^0$.
Let $j = \lceil \frac{k}{2} \rceil -1$.
We suppose that $\norm^0({\bf c}_1) \leq j$ and
$\norm^0({\bf c}_2) \leq j$ and derive a contradiction:
For $\ell = 1,2$ let branches $B_0^\ell, \dots, B_{j-1}^\ell$ and sets
$a^\ell \subseteq n_{3,i}$ exemplify this.

Let $a = a^1\cup a^2$ and let $\eta \in \rge(\val({\bf c}))$ be such that
for all $x \in (\dom(\eta) \cap \bigcup_{\ell=1,2}\bigcup_{i=0}^{j-1}
B^\ell_i)\setminus \dom(\eta({\bf c}))$ we have $\eta(x) \not\in a$.
But then for that $\ell \in \{1,2\}$ for which
$\eta \in \rge(\val({\bf c}_\ell))$ we get a contradiction.
Hence (1) follows for $\norm^1$
 $\norm^\ast$ increases or stays when taking subsets of $\val({\bf c})$,
and hence we have the analogous result for $\norm^{\frac{1}{2}}$.

Since the $k$-components of the creatures coincide, 
Part (2) follows from the behaviour of $\norm^{\frac{1}{2}}$ that
was shown in part (1) and from the requirements on $f$ in
Definition~\ref{1.7}(4): $ f(\frac{n}{2},k)
\geq f(n,k) -1$.
 \proofend

\section{Forcing with tree-creatures}\label{S2}

Now we define a notion of forcing with $\omega$-trees
$\langle {\bf c}_t \such t \in (T,\triangleleft_T) \rangle$  as conditions.
The nodes $t$ of these trees $(T,\triangleleft_T) = (\dom(p) , \triangleleft_p)$ 
and their immediate successors are described by certain
creatures ${\bf c}_t$   from Definition~\ref{1.5}. 
\nothing{The $\leq$-relation is quite complicated, but allows us to prove
many useful properties.
}
\smallskip

First we collect some general notation about trees.
The trees here are not the Aronszajn trees of the first section, but
trees $T$ of finite partial specialization functions, ordered by
$\triangleleft_T$ which is a subrelation of
$\subset$. Some of these trees will serve as forcing conditions.
\nothing{
The fact that
 the domains of these finite partial specialization functions are
not natural numbers, makes that \ref{1.3} does not describe them
well, and that we have to be careful with citations from
\cite{RoSh:470}. But still Sections 1.3 and 2.3 of that work
are useful, because they give the vocabulary about the
tree creatures and a information on sufficient conditions for properness.
The reasons here will be different, and especially \ref{3.2} will be new.
}

\begin{definition}\label{2.1}
\begin{myrules}
\item[(1)] A tree $(T,\triangleleft_T)$ is a set $T\subseteq \spec$,
such that for any $\eta \in T$, $(\{\nu \such \nu \triangleleft_T \eta\},
\triangleleft_T)$
is a finite linear order and such that
in $T$ there is  one least element, called the root, $\rt(T)$.
If $\eta \triangleleft_T \nu$ then $\eta \subset \nu$.
Every $\eta \in T\setminus \rt(T)$ can has just one
immediate 
$\triangleleft_T$-predecessor in $T$.
We shall only work with finitely branching trees.

\item[(2)]
We define the successors of $\eta$ in $T$, the restriction of 
$T$ to $\eta$, the splitting points
of $T$ and the maximal points of $T$ by

$$\suc_T(\eta) = \{ \nu \in T \such \eta \triangleleft_T \nu \wedge
\neg (\exists \rho \in T) (\eta \triangleleft_T \rho \triangleleft_T \nu) \},$$
$$T^{\langle \eta\rangle} = \{ \nu \in T \such \eta \trianglelefteq_T \nu \}.$$
$$\splitt(T)= \{ \eta \in T \such |\suc_T(\eta)| \geq 2 \},$$
$$\max(T)= \{ \nu \in T \such 
\neg (\exists \rho \in T)(\nu \triangleleft_T \rho\}.$$

\item[(3)] 
The $n$-th level of $T$ is
$$T^{[n]} = \{ \eta \in T \such \eta \mbox{ has $n$ 
$\triangleleft_T$-predecessors}\}.$$
The set of all branches through $T$ is
\begin{equation*}
\begin{split} \lim(T)=
\{ \langle \eta_k \such k <\ell \rangle \such & \ell \leq \omega 
\wedge (\forall k < \ell)( \eta_k \in T^{[k]})\\
&\wedge 
(\forall k <\ell-1)( \eta_k \triangleleft_T \eta_{k+1})\\
&\wedge \neg (\exists \eta_\ell \in T)(\forall k < \ell)( \eta_k 
\triangleleft_T \eta_\ell)\}.
\end{split}
\end{equation*}

A tree is well-founded if there are no  infinite branches through it.

\item[(4)]
A subset $F$ of $T$ is called a 
front of $T$ if every branch of $T$ passes through this set,
and the set consists of $\triangleleft_T$-incomparable elements.
\end{myrules}
\end{definition}

\begin{definition}\label{2.2}
We define a notion of forcing $Q= Q_{T_A}$. $p \in Q$ iff
\begin{myrules}

\item[(i)] $p$ is a function from a subset of $\spec=\spec^{T_A}$
(see Definition~\ref{1.12}) 
to $\omega$.

\item[(ii)] $p^{[]} = (\dom(p), \vartriangleleft_p)$ 
is a tree with $\omega$ levels,
the $\ell$-th level of which is denoted by $p^{[\ell]}$.

\item[(iii)] $p^{[]}$ has a root, the unique element of level 0,
called $\rt(p)$.

\item[(iv)] We let $$i(p) \stackrel{def}{=} 
\min\{ i \such |\dom(\rt(p))| \leq n_{2,i-1}\}.$$
Then for any $\ell < \omega$ 
and $\eta \in p^{[\ell ]}$ the set
$$\suc_p(\eta) = \{ \nu \in p^{[\ell +1]} \such \eta \subseteq \nu \}$$
is $\rge(\val({\bf c}))$ for a simple $(i(p) + \ell)$-creature ${\bf c}$ with 
base $\eta$. \nothing{
(Note that $n_{2,i} + n_{3,i} < n_{3,i+1}$.)} 
We denote this simple creature by ${\bf c}_{p,\eta}$ and let
${\bf c}^+_{p,\eta} = ({\bf c}_{p,\eta}, p(\eta))$.
Furthermore, we require
$p(\eta) \leq \norm^{\frac{1}{2}}({\bf c}_{p,\eta})$.

\item[(v)] If $\eta \in \dom(p)$ and $\nu \in \dom(p)$ and if
$\eta \cup \nu\in \spec$, then $\eta \cup \nu \in \dom(p)$.
It is a superset of both $\tau$ and of $\nu$, but in $\triangleleft_p$
it has only one predecessor.
Every $\eta \in \spec$ appears at most once in $T$.

\nothing{{\sf I added this item, because the work with the trees $T$ did not
go through smoothly with the $\subset$-relation.
So we work with an abstract $\triangleleft_T$ which has to fulfil some
requirements in order to make the projection functions (see below, 
in \ref{2.1}) unique.}}

\item[(vi)] For some $k < \omega$ for every $\eta \in p^{[k]}$ there is 
$\alpha < \omega_1$ and a finite $u \in T_A\setminus 
(T_A)_{<\alpha}$ such that for every 
$\omega$-branch $\langle \eta_\ell \such \ell < \omega \rangle$ 
of $p^{[]}$ satisfying $\eta_k = \eta$ we have 
$\bigcup_{\ell \in \omega} \dom(\eta_\ell) \setminus
u = (T_A)_{< \alpha}$.

\item[(vii)] For every $\omega$-branch
$\langle \eta_\ell \such \ell \in \omega \rangle$ of $p^{[]}$ we have
$\lim_{\ell \to \omega} \norm({\bf c}^+_{p,\eta_\ell}) =\omega$.
\end{myrules}

The order $\leq = \leq_Q$ is given by letting $p \leq q$ ($q$ is
stronger than $p$, we follow the Jerusalem convention) iff
$i(p) \leq i(q)$ and there is a projection $\pr_{p,q}$ which satisfies
\begin{myrules}
\item[(a)] $\pr_{q,p}$ is a function from $\dom(q)$ to $\dom(p)$.
\item[(b)] $\eta \in q^{[\ell]}
\Rightarrow \pr_{q,p}(\eta) \in p^{[\ell+i(q)-i(p)]}$.

\item[(c)] If $\eta_1 , \eta_2$ are both in $q^{[]}$
and if $\eta_1 \trianglelefteq_q \eta_2$,
then $\pr_{q,p}(\eta_1) \trianglelefteq_p \pr_{p,q}(\eta_2)$.
\item[(d)] $q(\eta) \geq p(\pr_{q,p}(\eta))$.
\nothing{{\sf We reversed the sign}}
\item[(e)] If $\eta \in q^{[]}$ then $\eta \supseteq \pr_{q,p}(\eta)$.
\nothing{{\sf I reversed the inclusion sign.}}
\item[(f)] If $\nu \in q^{[\ell]}$ and $\rho \in q^{[\ell+1]}$ and $
\nu \subseteq \rho$, $\pr_{q,p}(\nu)= \eta$, $\pr_{q,p}(\rho)= \tau$,
then $\dom(\tau) \cap \dom(\nu) = \dom(\eta)$.

\nothing{if $\nu \in q^{[\ell]}$ and $\rho \in q^{[\ell+1]}$ and $
\nu \subseteq \rho$, $\pr_{q,p}(\nu)= \eta$ then $\pr_{q,p}(\rho)=
\bigcup\{ \tau \such \tau \in \val(\eta) \wedge \rho \subseteq \tau\}$.}
\end{myrules}
\end{definition}

\nothing{
\renewcommand{\baselinestretch}{1} 
\small
\begin{figure}
  \begin{center}
    \ForceWidth{\textwidth}
    \BoxedEPSF{trees604A.eps} 
  \end{center}
\end{figure}
\renewcommand{\baselinestretch}{1.2}
\normalsize
}
\begin{definition}\label{2.3}
For $p \in Q$ and $\eta \in \dom(p)$ we let
$$p^{\langle \eta \rangle} = p \restriction \{ \rho \in \dom(p)
\such \eta \subseteq \rho \}.$$
\nothing{\item[(2)] 
For $p \in Q$, the set $\{ \eta_1,\dots,\eta_n \}$ is called a 
front of $\dom(p)$ if every branch of $p$ passes through this set,
and the set consists of $\subseteq$-incomparable elements of $\dom(p)$.
\end{myrules}}
\end{definition}

Let us give some informal description of the $\leq$-relation in $Q$:
The stronger conditions' domain is via $\pr_{q,p}$ mapped homomorphically w.r.t.\
the tree orders into $p^{\langle \pr_{q,p}(\rt(q))\rangle}$.
The stem can grow as well.  According to (b), the projection
preserves the levels in the trees but for
one jump in heights (the $\ell$'s in $p^{[\ell]}$), 
due to a possible lengthening of the
 stem. The partial specialization functions
sitting on the nodes of the tree are extended (possibly by more than one extension
to one function) in $q$ as to compared with the ones attached to the
image under $\pr$, but by (b)
the extensions
are so small and so few that it preserves the kind $i$ of the 
creature given by the node and its successors,  and according to (f)
the new part of the domain of the 
extension is disjoint from the domains of the
old partial specification
functions living higher up in the new tree.

\smallskip

Let us compare our setting with the forcings given in the book
\cite{RoSh:470}:
There the $\leq$-relation is based on a sub-composition
function (whose definition is
not used here, because we just deal with one
particular forcing notion)
 whose inputs are well-founded subtrees of the
weaker condition.
This well-foundedness condition  
\cite[1.1.3]{RoSh:470}  is not fulfilled: if we look at (e) and (f) in
the definition of $\leq$ we see that we have to look at all the branches
of $p$ that are in the range of $\pr_{q,p}$
in order to see whether some $\nu \in q^{[\ell]}$ fulfils (f) of the
definition of
$p \leq q$.
On the other hand, the projections shift all the levels 
by the same amount 
$i(q) - i(p)$, and are not arbitrary finite contractions as
in most of the forcings in the book \cite{RoSh:470}.

\nothing{
Nevertheless we present now parts framework from 
\cite[Section 1.3]{RoSh470} here,
because we still can use the nomenclature and 
many parts of proofs from there. We will add some new parts.
We try to write in a self-contained manner, so that the
not so interested reader need not look at the book.
}
\nc{\bH}{{\bf H}}
\nc{\TCR}{\mbox{TCR}}
\nc{\WCR}{\mbox{WCR}}
\renewcommand{\S}{{\mathcal S}}\nc{\C}{{\mathcal C}}
\nc{\rng}{\rge}
\renewcommand{\root}{\rt}
\newcommand{\nor}{\norm}
\renewcommand{\tree}{{\rm tree}}
\nc{\q}{{\mathbb Q}}
\nc{\rest}{\restriction}
\nc{\forces}{\Vdash}
\nc{\Wtil}{\name{W}}

\nothing{
\begin{definition}
\label{treecreature}(\ref[1.3.3]{RoSh:470}) 
\begin{enumerate}
\item A weak creature $t\in\WCR[\bH]$ is a {\em tree--creature} if $\dom(\val[
t])$ is a singleton $\{\eta\}$ and no two distinct elements of $\rng(\val[t])$
are $\vartriangleleft$--comparable;  
$\TCR[\bH]$ is the family of all tree--creatures for $\bH$. 
\item $\TCR_\eta[\bH]=\{t\in\TCR[\bH]: \dom(\val[t])=\{\eta\}\}$.
\item A sub-composition operation $\Sigma$ on $K\subseteq\TCR[\bH]$ is {\em a
tree composition} (and then $(K,\Sigma)$ is called {\em a tree--creating pair}
(for $\bH$)) if: 
\begin{description}
\item[(a)] if $\S\in [K]^{\textstyle{\leq}\omega}$, $\Sigma(\S)\neq\emptyset$
then $\S=\{s_\nu:\nu\in\hat{T}\}$ for some well founded quasi tree
$T\subseteq\bigcup\limits_{n<\omega}\prod\limits_{i<n}\bH(i)$ and a system
$\langle s_\nu: \nu\in\hat{T}\rangle\subseteq K$ such that for each finite
sequence $\nu\in\hat{T}$ 
\[s_\nu\in\TCR_\nu[\bH]\quad\mbox{ and }\quad\rng(\val[s_\nu])=\suc_T(\nu),\]
and
\item[(b)] if $t\in\Sigma(s_\nu: \nu\in\hat{T})$ then
$t\in\TCR_{\root(T)}[\bH]$ and $\rng(\val[t])\subseteq\max(T)$.
\end{description}
If $\hat{T}=\{\root(T)\}$, $t=t_{\root(T)}\in\TCR_{\root(T)}[\bH]$ and
$\rng(\val[t])=\max(T)$ then we will write $\Sigma(t)$ instead of
$\Sigma(t_\nu:\nu\in\hat{T})$. 
\item A tree-composition $\Sigma$ on $K$ is {\em bounded} if for each $t\in
\Sigma(s_\nu:\nu\in\hat{T})$ we have
\[\nor[t]\leq\max\{\nor[s_\nu]: (\exists \eta\in\rng(\val[t]))(\nu
\vartriangleleft\eta)\}.\] 
\end{enumerate}
\end{definition}
\begin{remark}
\label{treerem}
{\em
1)\ \ \ Note that sets of tree creatures relevant for tree compositions have a
natural structure: we identify here $\S$ with $\{s_{\nu(s)}: s\in\S\}$ where
$\nu(s)$ is such that $s\in\TCR_{\nu(s)}$ and $s_{\nu(s)}=s$.\\
2)\ \ \ To check consistency of our notation for tree creatures with that of
\ref{maindef} note that in \ref{treecreature}(3), if $s_\nu\in
\Sigma(s_{\nu,\eta}:\eta\in \hat{T}_\nu)$ for each $\nu\in\hat{T}$, $T$ is a
well founded quasi tree as in (3){\bf (a)} of \ref{treecreature} then
$T^*\stackrel{\rm def}{=}\bigcup\limits_{\nu\in\hat{T}} T_\nu$ is a well
founded quasi tree, $\hat{T}^*=\bigcup\limits_{\nu\in\hat{T}} \hat{T}_\nu$ and
$\langle s_{\nu,\eta}: \nu\in\hat{T}, \eta\in\hat{T}_\nu\rangle$ is a system
for which $\Sigma$ may be non-empty, i.e. it satisfies the requirements of
\ref{treecreature}(3){\bf (a)}.\\
3)\ \ \ Note that if $(K,\Sigma)$ is a tree--creating pair for $\bH$, $t\in
\TCR_\eta[\bH]$ then $\basis(t)=\{\eta\}$ and $\pos(\eta,t)=\rng(\val[t])$ (see
\ref{basis}). For this reason we will write $\pos(t)$ for $\pos(\eta,t)$ and
$\rng(\val[t])$ in the context of tree--creating pairs.\\ 
4)\ \ \ Tree--creating pairs have the properties corresponding to the niceness
and smoothness of creating pairs (see \ref{niceandsmo}, compare with
\ref{treesmooth}). 
}
\end{remark}
When dealing with tree--creating pairs it seems to be more natural to consider
both very special norm conditions and some restrictions on conditions of the
forcing notions we consider. The second is not very serious: the forcing
notions $\q^{\tree}_e(K,\Sigma)$ (for $e<5$) introduced in \ref{treeforcing}
below are dense subsets of the general forcing notions
$\q_{\C(\nor)}(K,\Sigma)$ (for suitable conditions $\C(\nor)$). We write the
definition of $\q^{\tree}_e(K,\Sigma)$ fully, not referring the reader to
\ref{maindef}, to show explicitly the way tree creating pairs work.
\begin{definition} 
\label{treeforcing} 
Let $(K,\Sigma)$ be a tree--creating pair for $\bH$.
\begin{enumerate}
\item We define the forcing notion $\q^{\tree}_1(K,\Sigma)$ by letting:
\noindent {\bf conditions } are sequences $p=\langle t_\eta: \eta\in T\rangle$
such that 
\begin{description}
\item[(a)] $T\subseteq\bigcup\limits_{n\in\omega}\prod\limits_{i<n}\bH(i)$ is
a non-empty quasi tree with $\max(T)=\emptyset$,
\item[(b)] $t_\eta\in\TCR_\eta[\bH]\cap K$ and $\pos(t_\eta)=\suc_T(\eta)$
(see \ref{treerem}(3)),
\item[(c)${}_1$] for every $\eta\in\lim(T)$ we have:
\[\mbox{the sequence }\langle\nor[t_{\eta \rest k}]:k<\omega, \eta\rest k\in
T\rangle\mbox{ diverges to infinity;}\]
\end{description}
\noindent{\bf the order} is given by:
\noindent $\langle t^1_\eta: \eta\in T^1\rangle\leq\langle t^2_\eta: \eta\in
T^2\rangle$ if and only if 
\noindent $T^2\subseteq T^1$ and for each $\eta\in T^2$ there is a well
founded quasi tree $T_{0,\eta}\subseteq (T^1)^{[\eta]}$ such that
$t^2_\eta\in\Sigma (t^1_\nu: \nu\in \hat{T}_{0,\eta})$.
If $p=\langle t_\eta:\eta\in T\rangle$ then we write $\root(p)=\root(T)$,
$T^p= T$, $t^p_\eta = t_\eta$ etc.
\item Similarly we define forcing notions $\q^{\tree}_e(K,\Sigma)$ for
$e=0,2,3,4$ replacing the condition {\bf (c)${}_1$} by {\bf
(c)${}_e$} respectively, where:
\begin{description}
\item[(c)${}_0$]  for every $\eta\in\lim(T)$:
\[\lim\sup\langle\nor[t_{\eta \rest k}]:k<\omega, \eta\rest k\in T\rangle
=\infty,\] 
\item[(c)${}_2$]  for every $\eta\in T$ and $n<\omega$ there is $\nu$ such
that  $\eta\vartriangleleft\nu\in T$ and $\nor[t_\nu]\geq n$,
\item[(c)${}_3$] for every $\eta\in T$ and $n<\omega$ there is $\nu$ such
that $\eta\vartriangleleft\nu\in T$ and 
\[(\forall\rho\in T)(\nu\vartriangleleft\rho\ \ \Rightarrow\ \
\nor[t_\rho]\geq n),\]  
\item[(c)${}_4$] for every $n<\omega$, the set 
\[\{\nu\in T: (\forall \rho\in T)(\nu\vartriangleleft\rho\quad\Rightarrow
\quad\nor[t_\rho]\geq n)\}\]
contains a front of the quasi tree $T$.
\end{description}
\item If $p\in\q^{\tree}_e(K,\Sigma)$ then we let $p^{[\eta]}=\langle t^p_\nu:
\nu\in (T^p)^{[\eta]}\rangle$ for $\eta\in T^p$.
\item For the sake of notational convenience we define partial order
$\q^{\tree}_{\emptyset}(K,\Sigma)$ in the same manner as
$\q^{\tree}_e(K,\Sigma)$ above but we omit the requirement {\bf (c)$_e$} (like
in \ref{conditions}; so this is essentially $\q_\emptyset(K,\Sigma)$).
\end{enumerate}
\end{definition}
\begin{remark}
{\em 
1)\ \ \ In the definition above we do not follow exactly the notation of
\ref{maindef}: we omit the first part $w^p$ of a condition $p$ as it can be
clearly read from the rest of the condition. Of course the missing item is
$\root(p)$. In this new notation the name $\dot{W}$ of \ref{thereal} may be
defined by 
\[\forces_{\q^{\tree}_e(K,\Sigma)}\Wtil=\bigcup\{\root(p): p\in
\Gamma_{\q^{\tree}_e(K,\Sigma)}\}.\] 
2)\ \ \ Note that 
\[\begin{array}{rr}
\q^{\tree}_4(K,\Sigma)\subseteq\q^{\tree}_1(K,\Sigma)\subseteq \q^{\tree}_0
(K,\Sigma)\subseteq\q^{\tree}_2(K,\Sigma)&\quad\mbox{ and}\\
\q^{\tree}_1(K,\Sigma)\subseteq\q^{\tree}_3(K,\Sigma)\subseteq
\q^{\tree}_2(K,\Sigma)\\ 
  \end{array}\]
but in general these inclusions do not mean ``complete suborders''. If the
tree--creating pair is t-omittory (see \ref{tomit}) then $\q^{\tree}_4(K,
\Sigma)$ is dense in $\q^{\tree}_2(K,\Sigma)$ and thus all these forcing
notions are equivalent. If $(K,\Sigma)$ is $\bar{2}$--big (see \ref{kbig})
then $\q^{\tree}_4(K,\Sigma)$ is a dense in $\q^{\tree}_1(K,\Sigma)$ (see
\ref{fronor2}). 
}
\end{remark}
\begin{definition}
\label{thick}
Suppose $(K,\Sigma)$ is a tree creating pair, $e<5$, $p\in\q^{\tree}_e(K,
\Sigma)$. A set $A\subseteq T^p$ is called {\em an $e$-thick antichain} (or
just {\em a thick antichain}) if it is an antichain in
$(T^p,\vartriangleleft)$ and for every condition $q\in\q^{\tree}_e(K,\Sigma)$
stronger than $p$ the intersection $A\cap\dcl(T^q)$ is non-empty.
\end{definition}
\begin{proposition}
\label{frontetc}
Suppose that $(K,\Sigma)$ is a tree--creating pair for $\bH$, $e<5$,
$p\in\q^{\tree}_e(K,\Sigma)$ and $\eta\in T^p$. Then:
\begin{enumerate}
\item $\q_e^{\tree}(K,\Sigma)$ is a partial order.
\item Each $e$-thick antichain in $T^p$ is a maximal antichain. Every front of
$T^p$ is an $e$-thick antichain in $T^p$.
\item If $e\in\{1,3,4\}$, $n<\omega$ then the set 
\[\begin{array}{ll}
B_n(p)\stackrel{\rm def}{=}\{\eta\in T^p: & \mbox{(i) } \ (\forall \nu\in
T^p)(\eta\trianglelefteq\nu\ \ \Rightarrow\ \ \nor[t_\nu]> n)\mbox{ but }\\
\ & \mbox{(ii) } \mbox{ no }\eta'\vartriangleleft\eta, \eta'\in T^p \mbox{
satisfies (i)}\}\\ 
\end{array}\]
is a maximal $\vartriangleleft$-antichain in $T^p$. If $e=4$ then $B_n(p)$ is
a front of $T^p$. 
\item For every $m,n<\omega$ the set
\[\begin{array}{ll}
F_n^m(p)\stackrel{\rm def}{=}\{\eta\in T^p: & \mbox{(i) } \ \nor[t_\eta]>n
\mbox{ and} \\   
\ & \mbox{(ii)}\ |\{\eta'\in T^p: \eta'\vartriangleleft\eta\ \&\
\nor[t_{\eta'}]>n\}|=m\}\\ 
\end{array}\]
is a maximal $\vartriangleleft$-antichain of $T^p$. If $e\in\{0,1,4\}$ then
$F^m_n(p)$ is a front of $T^p$.  
\item If $K$ is finitary (so $|\val[t]|<\omega$ for $t\in K$, see
\ref{morecreat}) then every front of $T^p$ is a front of $\dcl(T^p)$ and hence
it is finite.  
\item If $\Sigma$ is bounded then each $F^m_n(p)$ is a thick antichain of
$T^p$.  
\item $p\leq p^{[\eta]}\in \q_e^{\tree}(K,\Sigma)$\quad and\quad $\root(p^{[
\eta]})=\eta$. \QED
\end{enumerate}
\end{proposition}
\begin{remark}
\label{treesmooth}
{\em
One of the useful properties of tree--creating pairs and forcing notions
$\q^{\tree}_e(K,\Sigma)$ is the following:
\begin{description}
\item[$(*)_{\ref{treesmooth}}$] {\em Suppose that $p,q\in\q^{\tree}_e(K,
\Sigma)$, $p\leq q$ (so in particular $T^q\subseteq T^p$), $\eta\in T^q$ and
$\nu\vartriangleleft\eta$, $\nu\in T^p$.
Then $p^{[\nu]}\leq q^{[\eta]}$.}
\end{description}
}
\end{remark}
\begin{definition}
\label{AxA}
Let $p,q\in\q^{\tree}_e(K,\Sigma)$, $e<3$ (and $(K,\Sigma)$ a tree--creating
pair). We define relations $\leq^e_n$ for $n\in\omega$ by:
\begin{enumerate}
\item If $e\in\{0,2\}$ then $p\leq^e_{n} q$ {\em if and only if}:
\noindent $p\leq^e_0 q$ (in $\q^{\tree}_e(K,\Sigma)$) if $p\leq q$ and
$\root(p)=\root(q)$,
\noindent $p\leq^e_{n+1} q$ (in $\q^{\tree}_e(K,\Sigma)$) if $p\leq^e_0 q$ and
if $\eta\in F^0_n(p)$ (see \ref{frontetc}(4)) and $\nu\in T^p$,
$\nu\trianglelefteq\eta$ then $\nu\in T^q$ and $t^q_\nu=t^p_\nu$. 
\item The relations $\leq^1_{n}$ (on $\q^{\tree}_1(K,\Sigma)$) are defined by:
\noindent  $p\leq^1_0 q$ (in $\q^{\tree}_1(K,\Sigma)$) if and only if $p\leq
q$ and $\root(p)=\root(q)$,
\noindent $p\leq_{n+1}^1 q$ (in $\q^{\tree}_1(K,\Sigma)$) if $p\leq^1_0 q$ and
if $\eta\in F^0_n(p)$ (see \ref{frontetc}(4)) and $\nu\in T^p$,
$\nu\trianglelefteq\eta$ then $\nu\in T^q$, $t^p_\nu=t^q_\nu$, and
\[\{t^q_\eta:\eta\in T^q\ \&\ \nor[t^q_\eta]\leq n\}\subseteq\{t^p_\eta:\eta
\in T^p\}.\]  
\item We may omit the superscript $e$ in $\leq^e_n$ if it is clear in which
of the forcing notions $\q^{\tree}_e(K,\Sigma)$ we are working.
\end{enumerate}
\end{definition}
\begin{proposition}
\label{fusAxA}
Let $(K,\Sigma)$ be a tree--creating pair for $\bH$, $e<3$.
\begin{enumerate}
\item The relations $\leq^e_n$ are partial orders on $\q^{\tree}_e(K,\Sigma)$
stronger than $\leq$. The partial order $\leq^e_{n+1}$ is stronger than
$\leq^e_n$. 
\item Suppose that conditions $p_n\in\q^{\tree}_e(K,\Sigma)$ are such that
$p_n\leq_{n+1}^e p_{n+1}$. 
Then $\lim\limits_{n\in\omega}p_n\in\q^{\tree}_e(K,\Sigma)$ and $(\forall 
n\in\omega)(p_n\leq^e_{n+1}\lim\limits_{n\in\omega} p_n)$ (where the limit
condition $p=\lim\limits_{n\in\omega}p_n$ is defined naturally; $T^p
=\bigcap\limits_{n\in\omega} T^{p_n}$). \QED 
\end{enumerate}
\end{proposition}
}

\begin{definition}\label{2.4}
\begin{myrules}
\item[(1)] $p \in Q$ is called normal iff for every $\omega$-branch
$\langle \eta_\ell \such \ell \in \omega \rangle$ of $p^{[]}$ the sequence
$\langle \norm({\bf c}^+_{p,\eta_\ell}) \such \ell \in \omega \rangle$ is
non-decreasing.
\item[(2)] $p \in Q$ is called smooth iff in clause (vi) of 
Definition~\ref{2.2} the number $k$ is 0 and $u$ is empty.

\item[(3)] $p \in Q$ is called weakly smooth iff in clause (vi) of 
Definition~\ref{2.2} the number $k$ is 0.
\end{myrules}
\end{definition}

\begin{remark}\label{2.5} If $p \in Q$ is smooth then there is some 
$\alpha < \omega_1$
such that for every $\omega$-branch $\langle \eta_\ell \such \ell \in \omega 
\rangle$ of $p^{[]}$ we have $\bigcup_{\ell < \omega} \dom(\eta_\ell)
= T_{<\alpha}$. This $\alpha$ is denoted by $\alpha(p)$.
 \end{remark}

\begin{fact}\label{2.6}
\begin{myrules}
\item[(1)]
Suppose that we have strengthened \ref{2.2}(f) to: 
 If $\nu \in q^{[\ell]}$ and $\rho \in q^{[\ell+1]}$ and $
\nu \subseteq \rho$, $\pr_{q,p}(\nu)= \eta$, $\tau \supseteq \eta$,
$\tau \in p$,
then $\dom(\tau) \cap \dom(\nu) = \dom(\eta)$.
So, 
we replaced  $\pr_{q,p}(\rho)= \tau$ by 
$\tau \supseteq \eta$, $\tau \in p$ and thus have information on the
$\tau \in \dom(p)\setminus \rge(\pr_{q,p})$.
Then:
If $p \leq q$ and $\eta \in \dom(p)$, 
$\nu \in \dom(q)$, $\eta = \pr_{q,p}(\nu)$ and 
$\eta \vartriangleleft \tau \in \dom(p)$, then
$\dom(\nu) \cap \dom(\tau) = \dom(\eta)$.
\item[(2)]
If $p \leq q$ and $p$ is weakly smooth then\\
 $\nu \in \dom(q) \rightarrow
\dom(\nu) \cap (T_{<\alpha(p)}\cup u)
=\dom(\pr_{q,p}(\nu))$.
\end{myrules}
\end{fact}
\proof (1) follows from  clauses (e) and the stronger form of clause (f)
of the definition
of $p \leq_Q q$. (2): If $p$ is weakly smooth, \ref{2.2}(f) and its stronger
form from the premise of (1) coincide, and hence (2) follows from (1),
because each branch of $p$ has the same union of domains.
\proofend

\begin{definition}\label{2.7}
For $0 \leq n < \omega$ we define 
the partial order $\leq_n$ on $Q$ by letting $p \leq_n q$ iff
\begin{myrules}
\item[(i)] $p \leq q$,
\item[(ii)]
$i(p) = i(q)$,
\item[(iii)] $p^{[\ell]} = q^{[\ell]}$ for $\ell \leq n$, 
and $p \restriction \bigcup_{\ell \leq n} p^{[\ell]} =
q \restriction \bigcup_{\ell < n} q^{[\ell]}$, in particular $\rt(p)=\rt(q)$,
\item[(iv)] if $\pr_{q,p}(\eta) = \nu$, then 
\begin{myrules}
\item[--] $\eta= \nu$ and ${\bf c}^+_{q,\eta} =
{\bf c}^+_{p,\nu}$
\item[--] or $\norm({\bf c}^+_{q,\eta}) \geq n$.
\end{myrules}
\end{myrules}
\nothing{For $n=0$ we let $\leq_0 = \leq$.}
\end{definition}

We state and prove some basic properties of the notions defined above.

\begin{claim}\label{2.8}

\begin{myrules}
\item[(1)]If $p \leq q$, then $\pr_{q,p}$ is unique. 
\item[(2)] If $p \in Q$ and $\ell \in \omega$ then $|p^{[\ell]}|
< n_{1,i(p) + \ell}$.
\item[(3)] $(Q,\leq_Q)$ is a partial order.

\item[(4)] If $p \leq q$ and $\pr_{q,p}(\eta) = \nu$, 
then $i({\bf c}_{q,\eta})=
i({\bf c}_{p,\nu})$.
\item[(5)] If $p \leq q$ and $\pr_{q,p}(\eta) = \nu$, 
then $\norm^0({\bf c}_{q,\eta})\leq \norm^0({\bf c}_{p,\nu})$.

\item[(6)] $(Q,\leq_n)$ is a partial order.
\item[(7)] $p \leq_{n+1}q \rightarrow p \leq_n q \rightarrow p\leq q$.
\item[(8)] If ${\bf c}$ is a simple $i$-creature with 
$k \leq \norm^{0}({\bf c})$, then there 
is a simple $i$-creature ${\bf c}'$ with
$k=\norm^0({\bf c}')$ and $\val({\bf c}') \subseteq \val({\bf c})$.
\item[(9)] For every $p \in Q$ there is a $q \geq p$ such that for
all $\eta$ and $\nu$ 
$$\pr_{q,p}(\eta) = \nu \rightarrow \norm^0({\bf c}_{q,\eta})
=\min\{\norm^0({\bf c}_{p,\rho}) \such \nu \trianglelefteq
\rho \in p \}.$$

\item[(10)] For every (not necessarily normal)
$p$ we have that $\lim_{n \to \omega} \min\{\norm({\bf c}^+ _{p,\eta})
\such \eta \in p^{[n]}\} = \infty$.
 
\item[(11)] If $p \in Q$ and $\eta \in p^{[\ell]}$ then $|\dom(\eta)| <
n_{2,i(p) + \ell-1}$ or $\ell =0$ and $i(p)=0$ and
$\eta=\emptyset$.
\end{myrules}
\end{claim}

\proof
(1) By induction on $\ell$ we show that $\pr_{q,p} \restriction
\bigcup_{\ell' \leq \ell} p^{[\ell']}$ is unique: It is easy to see
that $\pr_{q,p}(\rt(q))$ is the $\subseteq$-maximal element of $p$ that is
a subfunction of $\rt(q)$. By Definition~\ref{2.2}(5) such a maximum exists.
Then we proceed level by level in $q^{[]}$, and
again Definition~\ref{2.2}(5) yields uniqueness of $\pr_{q,p}$.

(2) This is also proved by induction on $\ell$.
Note that for $\eta \in p^{[\ell]}$ we have that $|\rge(\val((\eta))| \leq
n_{1,i(p) + \ell -1}$. We have $|p^{[0]}| =1$ and
by Definition~\ref{1.5}(c), 
$|p^{[\ell +1]}|\leq |p^{[\ell]}| \cdot n_{1,i(p) + \ell} \leq$
$n_{1,i(p) +\ell} \cdot n_{1,i(p)+\ell} \leq n_{1,i(p)+\ell+1}$,
by equation~\eqref{1.3}.

(3) Given $p \leq q$ and $q \leq r$ we define $\pr_{r,p}=
\pr_{q,p} \circ \pr_{r,q}$. It is easily seen that this function is as 
required.

(4) Let $\ell$ be such that $\eta \in q^{[\ell]}$. Then $i({\bf c}_{q,\eta}) 
= i(q) + \ell$ and $\nu \in p^{[\ell +i(q) - i(p)]}$.
Hence $i({\bf c}_{p,\nu}) = i(p) + \ell + i(q) -i(p) = i(q) + \ell$.

(5) Suppose $\norm^0({\bf c}_{q,\eta}) 
> \norm^0({\bf c}_{p,\nu})$.
Let $k = \norm^0(c_{q,\eta})$ and let $i=i({\bf c}_{q,\eta}) = 
i({\bf c}_{p,\nu})$. Suppose that $a \subseteq n_{3,i}$ and the branches $B_0,
\dots, B_{k-1}$ of $T$ exemplify that $\norm^0({\bf c}_{p,\nu}) < k$.
Hence for all $\tau \in \suc_p(\nu)$ 
\begin{myrules}
\item[($\alpha$)]
there is $x \in (\dom(\tau) 
\cap \bigcup_{\ell=0}^{k-1} B_\ell) \setminus \dom(\nu)$
such that $\tau(
x) \in a$, or
\item[$(\beta)$] 
$|\dom(\tau)| \geq \frac{n_{2,i}}{2^k}$.
\end{myrules}

Let $\tau \in \suc_q(\nu)$, and $\pr_{q,p}(\tau) = \tau'$.
Suppose $(\alpha$) is the case for $\tau'$. Then the same $a$ and $B_0, \dots, B_{k-1}$
exemplify $(\alpha)$ for $\tau$ and ${\bf c}_{q,\eta}$, because
 we have $\eta \supseteq \nu = \pr_{q,p}(\eta) $ and 
 $\pr_{q,p}[\suc_q(\eta)] \subseteq \suc_p(\nu)$.
The same $x$ will show that $\exists x \in (\dom(\tau) 
\cap \bigcup_{\ell=0}^{k-1} B_\ell) \setminus \dom(\eta)$
such that $\tau(
x) \in a$, if we verify that $x \not\in \dom(\eta)$. 
But we have for all $\tau \in \suc_q(\eta)$ that
$\dom(\eta) \cap \dom(\tau) = \dom(\nu)$ by \ref{2.2}(f), and hence
$x\not\in \dom(\eta)$.

Suppose $(\beta)$ is the case for $\tau'$. Then $\tau' \in \suc_p(\nu)$ and
$\tau \supseteq \tau'$, and hence
$|\dom(\tau)| \geq \frac{n_{2,i}}{2^k}$.

\smallskip

(6) Suppose that $p \leq_n q \leq_n r$ and $\pr_{r,q} (\sigma) 
= \eta$ and $\pr_{q,p}(\eta) = \nu$.
By (1) and (3) we have that $\pr_{r,p}(\sigma) = \nu$, 
and now it is easy to check the requirements for $p \leq_n r$.

\smallskip

(7) Obvious.

\smallskip

(8) We may assume that $\norm^0({\bf c})
> k$, because otherwise $\bf c$ itself is as required.
Look at
\begin{equation*}\begin{split}
Y = \{{\bf d} \such & {\bf d} \mbox{ is a simple $i$-creature and }
\val({\bf d}) \neq \emptyset \mbox{ and }\\&
\norm^0({\bf d}) \geq k \mbox{ and } \val({\bf d}) \subseteq \val({\bf c})\}.
\end{split} \end{equation*}
 
Since ${\bf c} \in Y$, it is non-empty, and it has a member $\bf d$ with 
a minimal number of elements. We assume towards a contradiction that 
$\norm^0({\bf d}) > k$.
We choose $\eta^* \in \rge(\val({\bf d}))$. We let $\rge(\val({\bf d}^*)) = 
\rge(\val({\bf d})) \setminus 
\{\eta^*\}$.

Claim: ${\bf d}^\ast \neq \emptyset$. Otherwise we choose
$x \in \dom(\eta^\ast) \setminus \dom(\eta({\bf c}))$, and such 
an $x$ exists by clause (d) of \ref{1.3}(1). Now we let $B_0$ 
be a branch of $T$ to which $x$ belongs and set 
$a = \{ \eta^\ast(x)\}$. They witness that 
$\norm^0({\bf c}) \not\geq 1$, so $\norm^0({\bf c}) =0$, 
which contradicts the
assumption that $\norm^0({\bf c}) = k > 0$.

Claim: $\norm^0({\bf d}^\ast) \geq k$. Otherwise there are branches $B_0,\dots,B_{k-1}$
and a set $a\subseteq n_{3,i}$ witnessing  $\norm^0({\bf d}^\ast) \not\geq k$. 
Let $x \in \dom(\eta^\ast) \setminus \dom(\eta({\bf c}))$ 
and let $B_k$ be a branch
such that $x \in B_k$ and set $a' = a \cup\{\eta^\ast(x)\}$.
The $B_0,\dots , B_k$ and $a'$ witness that $\norm^0({\bf c}) \not\geq k+1$.
Hence ${\bf d}^\ast$ is a member of $Y$ with fewer elements than ${\bf d}$, contradiction.


\smallskip

(9) Follows from (8). We can even take $q \subseteq p$.
First see:  For no $m$ such the set
 $\{\eta\in p$ such that for densely (in $p^{[]}$) many $\eta'
\trianglerighteq_{p} \eta$ we have that $\norm^0({\bf c}_{p,\eta'}) < m\}$.
is anywhere dense.
Otherwise we an choose a branch 
$\langle \eta_\ell \such \ell \in \omega \rangle$ 
such that there is some $m \in \omega$ such that
for all $\ell < \omega$, $\nor^0({\bf c}_{p,\eta_\ell})<m$.

Now we choose by induction of $\ell$, $\dom(q_\ell) \subseteq \dom(p)$, such
that $q_\ell$ is 
has no infinite branch and hence is finite, though
we do not have a bound on its height.

\smallskip
First step:
Say $\min\{ \norm^0({\bf c}_{p,\eta}) \such \eta \in \dom(p)\} = k$ and
it is reached in $\eta \in \dom(p)$.
We take $q^{[0]}=\{\eta\}$.

([1])
Then we take for any $\eta' \in \rge(\val({\bf c}_{p,\eta}))$ some $\eta'' 
\supseteq \eta'$ such 
that $\eta'' \in p$ and such that for all $\tilde{\eta} \supseteq
\eta''$, if $\tilde{\eta} \in p$ then $\norm^0({\bf c}_{p,\tilde{\eta}}) \geq 
k+1$. By the mentioned nowhere-density result, this is possible.
We put such an $\eta''$ in $q^{[\ell]}$, if it is in $p^{[\ell]}$.

([2])
Then we look at the $\nu$ in the
 branch between $\eta$ and $\eta''$ in $\dom(p)$.
If $\norm^0({\bf c}_{p,\nu}) > k$ we take according to (8) a subset of
$\rge(\val({\bf c}_{p,\nu}))$ with norm $k$ and put this into $q$.
We have to put successors to all $\nu' \in \rge(\val({\bf c}_{p,\nu}))$
for all $\nu$ in question into $q_\ell$.
This is done as in ([1]), applied to $\nu$ instead of $\eta$.
With all the $\nu$ in this subset we do the procedure in ([1]), 
and repeat and repeat it.  In 
finitely many (intermediate)
steps we  reach a subtree $\dom(q_\ell)$ of $\dom(p)$ without any
$\omega$-branches 
such that all its leaves fulfil
 $\eta'' \in p$ and such that for all $\tilde{\eta} \supseteq
\eta''$, if $\tilde{\eta} \in p$ then $\norm^0({\bf c}_{p,\tilde{\eta}}) \geq 
k+1$, and all its nodes $\eta$ fulfil $\norm^0({\bf c}_{q_1,\eta}) \geq k$.
By K\"onig's lemma, this tree $q_\ell$ is finite.

([3]) With the leaves of $q_\ell$ and $k+2$ instead of $k+1$, we repeat the choice procedure in ([1]) and 
([2]).
We do it successively for all $k\in \omega$. The union of the $q_\ell$, $\ell \in \omega$, is 
a $q$ as desired in (9).

\smallskip
(10) This follows from
 K\"onig's lemma: Since $p^{[]}$ is finitely branching, 
there is a branch though every infinite subset.

\smallskip

(11) Follows from Definitions~\ref{1.5} and \ref{2.2}.
\proofend

The next lemma states that $Q$ fulfils some fusion property:

\begin{lemma}\label{2.9}
\nothing{Let $f \colon \omega \to \omega$ be an increasing function such that
$\lim f(i) = \infty$. Let $q_i \in Q$, $i < \omega$.
We choose $n_i$ such that
$(\forall n \geq n_i)(\forall \eta \in q_i^{[n]}) \norm^0({\bf c}_{q,\eta})
\geq f(i)$.}
Let $\langle n_i \such i \in \omega \rangle$ be a strictly
increasing sequence of natural numbers.
 We assume that for every $i$,
$q_i \leq_{n_i} q_{i+1}$, and we set $n_{-1} = 0$.
Then $q=
\bigcup_{i < \omega} \bigcup_{n_{i-1} \leq n < n_{i}} q_i^{[n]} \in Q$ 
and for all $i$, $q \geq_{n_i} q_i$.
\end{lemma}

\proof Clear by the definitions.

The fusion lemma is usually applied in the following setting:

\begin{conclusion}\label{2.10}
Suppose $p \in Q$ us given and we are to find $q \geq p$
such that $q$ fulfils countably many tasks.
For this it is enough to find for any single task and any
$p_0$ and $k^\ast \in \omega$ some $q \geq_{k^\ast} p_0$ 
that fulfils the task.
\end{conclusion}

Now we want to 
fill up the domains of the partial specialization functions
and to show that smooth conditions are dense:

\begin{lemma}\label{2.11}
If $p \in Q$ and  
$m < \omega$ then for some smooth $q \in Q$ we have $p \leq_m q$. Moreover, if
$\bigcup \{ \dom(\eta) \such \eta \in p^{[]} \} \subseteq T_{<\alpha}$ then
we can demand that
$\bigcup \{ \dom(\eta) \such \eta \in q^{[]} \} = T_{<\alpha}$.
Moreover, $\eta \in q^{[]}$ implies
$\norm^1({\bf c}_{q,\eta}) \geq
\norm^1({\bf c}_{p,\pr_{q,p}(\eta)}) -1$, 
and $q(\eta)=p(\pr_{q,p}(\eta))$ implies that
$\norm({\bf c}^+_{q,\eta}) \geq \norm({\bf c}^+_{p,\pr_{q,p}(\eta)})-1$.
\end{lemma}

\proof We first use the definition of
$p \in Q$:
By item (v) there is
some $k < \omega$ for every $\eta \in p^{[k]}$ there is 
$\alpha(\eta) < \omega_1$ and a finite $u_\eta \in T_A\setminus 
(T_A)_{<\alpha(\eta)}$ such that for every 
$\omega$-branch $\langle \eta_\ell \such \ell < \omega \rangle$ 
of $p^{[]}$ satisfying $\eta_k = \eta$ we have 
$\bigcup_{\ell \in \omega} \dom(\eta_\ell) \setminus
u_\eta = (T_A)_{< \alpha(\eta)}$.
We fix such a $k$ and such $u_\eta$'s.

We can find $n$ such that
\begin{myrules}
\item[$(\ast)_1$] $m \leq n < \omega$, $k \leq n$,
\item[$(\ast)_2$]
 $|\bigcup_{\eta \in p^{[k]}}
{u_\eta}| <n$,
\nothing{
\item[$(\ast)_2$] for every $\eta \in p^{[n]}$ there is $\alpha(\eta)$ and a
finite $u_\eta \subseteq T \setminus T_{<\alpha(\eta)}$
such that
if $\langle \eta_k \such k \in \omega \rangle$ is a $\omega$-branch of
$p$ and $\eta = \eta_n$, then $\bigcup_{k \in \omega} \dom(\eta_k \setminus
u_\eta) = T_{<\alpha(\eta)}$.}
\item[$(\ast)_3$] for every $\nu \in p^{[n]}$, we have
$\norm^0({\bf c}_{p,\nu}) > m$,
\item[$(\ast)_4$] if $\eta \in p^{[n]}$, $\eta \subseteq \nu \in p$ then
$\dom(\nu) \setminus  \dom(\eta)$ is disjoint from $u$.
\end{myrules}

For each $\eta \in p^{[n]}$ let
$w^+_\eta =
\{ \nu \such \eta \vartriangleleft_T \nu \in \dom(p) \wedge
\norm^1({\bf c}_{p,\nu}) > 
\norm^1({\bf c}_{p,\eta})\}$
and $w_\eta = \{ \nu \in w_\eta^+ \such 
(\not\exists \rho)(\eta \trianglelefteq_T \rho \vartriangleleft_T \nu \wedge
\rho \in w^+_\eta \}$.
So $w := \bigcup \{ w_\eta \such \eta \in p^{[n]}\}$ is a front of $p^{[]}$.
For each $\nu \in w$ let $\nu \in p^{[\ell(\nu)]}$ (so $\ell(\nu) \geq n$)  and let $\alpha(\nu) 
= \alpha(\eta)$ and $u_\nu= u_\eta$ when $\nu \in w_\eta$.
Let $\{x^\nu_\ell \such \ell \in (\ell(\nu),\omega)\}$ 
enumerate $T_{<\alpha }
\setminus (T_{<\alpha(\nu)}\cup u_\nu)$ without repetition.
For each $\ell$, $\nu$,
\begin{equation}
\tag*{$\boxtimes$} 
|\{ \rho \in u_\nu \such x^\nu_\ell \triangleleft_T \rho 
\}| + 
|\{x^\nu_k \such \ell(\nu)<k<\ell \}| < \ell.
\end{equation}
We let
\begin{equation*}
\begin{split}
q_\ast^{[]} = \; \{ \rho \such & 
\mbox{ (a) for some } \nu \in w (\rho \subseteq
 \nu \vee \nu \subseteq \rho \in p^{[]},
\alpha= \alpha(\nu)\\
&  \mbox{ or (b) }
\nu \subseteq \rho, \; \alpha(\nu) < \alpha \mbox{ and for some }
\ell > \ell(\nu) \\
& \makebox{} \mbox{ and } \exists \tau \in p^{[\ell]} (\nu \subseteq \tau\\ 
&
\dom(\rho) =\dom(\tau) \cup \{x_k \such \ell(\nu) < k< \ell 
\mbox{ and } \rho(x^\nu_k) < n_{3,i(p) +k})\}
\end{split}
\end{equation*}
and choose $q^{[]} \subseteq q_\ast^{[]}$ 
by successively climbing upwards in the levels of $p^{[]}$, using first Claim~\ref{1.10}
for the immediate successors of an already chosen node, and then using Claim~\ref{1.11}
to make these new successors the bases of the creatures attached to them.
We choose $q$ rich enough  as in Claim~\ref{1.10} 
but also small enough as to have sufficienty high $\norm^\ast({\bf c}_{q,\rho})$.
We set
$q(\rho) = p(\tau)$ if $\rho,\tau$ are as above.

For checking the conditions for $p \leq q$ and on the norms   note that 
$\boxtimes$ above gives clause of the premises of Claim~\ref{1.10}
on a given level  and of Claim~\ref{1.11} on its successor level.
By the choice of $q$, it is smooth.

\begin{conclusion}\label{2.12}
Forcing with $Q$ specializes $T_A$.
\end{conclusion}

\section{Decisions taken by the tree creature forcing}
\label{S3}

In this section we prove that $Q$ is proper.
Indeed we prove that $Q$ has ``continuous reading of
names'' (this is the property stated in \ref{3.8}), which implies 
Axiom A (see \cite{Baumgartner}) and properness.

\begin{claim} \label{3.1}
\begin{myrules}
\item[(1)] If $p \in Q$ and $\{\eta_1,\dots,\eta_n\}$ is a front of $p$,
then $\{p^{\langle\eta_1\rangle},\dots,p^{\langle\eta_n\rangle}\}$ is
predense above $p$.
\item[(2)] If  $\{\eta_1,\dots,\eta_n\}$ is a front of $p$ and
 $p^{\langle\eta_\ell\rangle}
\leq q_\ell \in Q$ for each $\ell$, then there is $q \geq p$ with
$\{\eta_1,\dots,\eta_n\}\subseteq q^{[]}$ such that
for all $\ell$ we have that
$q^{\langle \eta_\ell\rangle} = q_\ell$.
Hence $\{q^{\langle \eta_\ell\rangle} \such 1 \leq \ell \leq n \}$ 
is predense above $q$.
\end{myrules}
\end{claim}

\begin{claim}\label{3.2}

If $p \in Q$
and 
$X \subseteq \dom(p)$ is upwards closed in $\triangleleft_p$,
and $\forall \eta \in \dom(p) \norm^0({\bf c}_{p,\eta}) > 0$, then there
is some $q$ such that
\begin{myrules}
\item[(a)] $ p \leq_0 q$, and either $(\exists \ell) q^{[\geq \ell]}
\subseteq X$ or $\dom(q) \cap X = \emptyset$,
\item[(b)] $\dom(q) \subseteq \dom(p)$ and $q = p \restriction \dom(q)$,
\item[(c)] for every $\nu \in \dom(q)$, if
${\bf c}_{q,\nu} \neq {\bf c}_{p,\nu}$, then 
$\norm^2({\bf c}_{q,\nu}) \geq \norm^2({\bf c}_{p,\nu})-1$
and $\norm({\bf c}^+_{q,\nu}) \geq \norm({\bf c}^+_{p,\nu})-1$.
\end{myrules}
\end{claim}

\proof 
We will choose $\dom(q) \subseteq \dom(p)$ and then let
$q = p \restriction \dom(q)$. 
For each $\ell$
we first 
choose by downward induction on $j \leq \ell$
subsets  $X_{\ell,j} \subseteq p^{[\leq \ell]}$
and a colouring $f_{\ell,j}$ of $X_{\ell,j}\cap p^{[j]}$ 
with two colours, 0 and 1.
The choice is performed in such a way that
$X_{\ell,j-1} \subseteq X_{\ell,j}$ and such that
$p^{[i]} \subseteq X_{\ell,j}$ for $i \leq j$.

\nothing{ looking first at the
highest level $p^{[\ell]}$. We choose $q_\ell$ in $\ell-1$ steps,
$q_{\ell,0} = p^{[\leq \ell]} \supseteq q_{\ell,1} \mdots
\supseteq q_{\ell, \ell-1} = q_\ell$.
Suppose that $q_{\ell,k} $ is chosen and we are to choose $
q_{\ell,k+1}$.}

We choose $X_{\ell,\ell} = p^{[\leq \ell]}$ and 
for $\nu \in p^{[\ell]}$ we set $f_{\ell,\ell}(\nu)
= 0$ iff $(\exists \ell')
(p^{\langle \nu \rangle})^{[\geq\ell']} \subseteq X$
and $f_{\ell,\ell}(\nu)=1$ otherwise.

Suppose that $X_{\ell,j}$ and $f_{\ell,j}$ are chosen.
For $\eta \in p^{[j-1]}\cap X_{\ell,j}$ we have 
\begin{equation*}
\begin{split}
\rge(\val({\bf c}_{\eta,p})) = &
\{ \nu \in \rge(\val({\bf c}_{\eta,p})) \such f_{\ell,j} =0\} \cup \\
& \{ \nu \in \rge(\val({\bf c}_{\eta,p})) \such f_{\ell,j} =1\} 
\end{split}
\end{equation*}
Note that the sets would be all the same if we intersect with
$X_{\ell,j}$, because $p^{[j]} \subseteq X_{\ell,j}$.
By Claim \ref{1.12} at least one of the two sets gives a creature
$\bf c$
with $\norm^2({\bf c}) > \norm^2({\bf c}_{\eta,p}) -1$.
So we keep in  $X_{\ell,j-1}\cap p^{[j]}$ 
only those of the majority colour
and close this set downwards in $p^{[]}$. This is $X_{\ell,j-1}$.
We colour the points on 
$p^{[j-1]}\cap X_{\ell,j-1}$ 
 with $f_{\ell,j-1}$
according to these majority colors, i.e.,
$f_{\ell, j-1}(\eta) = i$ iff
$\{ \nu \in \rge(\val({\bf c}_{\eta,p})) \such
f_{\ell,j}(\nu)  = i \} \subseteq X_{\ell,j-1}$.
We work downwards until we come to the root of $p$ and
keep $f_{\ell,0}(\rt(p))$ in our memory.

We repeat the procedure of the downwards
induction on $j$ for larger and larger $\ell$.

\smallskip

 If there is one $\ell$ where the root
got colour 0, we are, because $X$ is upwards closed,
 in the first case of the alternative in the conclusion 
(a). If for all $\ell$ the root got colour $0$, we have for all $\ell$ 
finite subtrees $t$ such that for all $\nu \in t$,
$p^{\langle \nu \rangle} \cap t$ has sufficiently high norm at its root.
By K\"onig's Lemma
 (initial segments of trees are taking from finitely 
many possibilities) we build a 
condition $q$ that all of its nodes are not in $X$, and thus (a) 
is proved. The item
(b) is clear. Item (c) follows from our choice of $q$ and from
 \ref{1.12}.
\proofend

The next claim is very similar to \ref{3.2}. We want to find
$q \geq_m p$, and therefore we have to weaken
the homogeneity property in item $(a)$ of \ref{3.2}.

\begin{claim}\label{3.3}

If $p \in Q$
$k^\ast \in \omega$, 
and 
$X \subseteq \dom(p)$ is upwards closed,
and $\forall \eta \in \dom(p) \norm^0({\bf c}_{p,\eta}) > 0$, then there
is some $q$ such that
\begin{myrules}
\item[(a)] $ p \leq_{k^\ast} q$, and there is a front
$\{\nu_0,\dots \nu_s\}$ of $p$ which is contained in $q$ and
whose being contained in $q$ ensures $p \leq_{k^\ast} q$,
and such that for all $\nu_i$ we have:
either $(\exists \ell) (q^{\langle \nu_i\rangle})^{[\geq \ell]}
\subseteq X$ or $\dom(q^{\langle \nu_i \rangle}) \cap X = \emptyset$,
\item[(b)] $\dom(q) \subseteq \dom(p)$ and $q = p \restriction \dom(q)$,
\item[(c)] for every $\nu \in \dom(q)$, if
${\bf c}_{q,\nu}^+ \neq {\bf c}_{p,\nu}^+$, then 
$\norm^2({\bf c}_{q,\nu}) \geq \norm^2({\bf c}_{p,\nu})-1$
and  $\norm({\bf c}^+_{q,\nu}) \geq \norm({\bf c}^+_{p,\nu})-1$.
\end{myrules}
\end{claim}

\proof We repeat the proof of \ref{3.2} for each $p^{\langle \nu_i\rangle}$.
\proofend

Now for the first time we make use of the coordinate
$k({\bf c}^+)$ of our creatures. The next lemma states that
the creatures have the halving property (compare
to \cite[2.2.7]{RoSh:470}).

\begin{definition}\label{3.4}

$Q$ has the halving property, iff there is a
function $\halv \colon K^+ \to K^+$ with the following properties:
\begin{myrules}
\item[(1)]
$\halv({\bf c}^+) = ({\bf c},k(\halv({\bf c}^+)))$,
\item[(2)]
$\norm(\halv({\bf c}^+))  \geq  
 \norm({\bf c}^+)-1$,
\item[(3)]
if ${\bf c}'$ is a simple creature and $ k\geq k(\halv({\bf c}^+))
$ and 
$\norm({\bf c}',k) > 0$, then $\norm({\bf c}', k({\bf c}^+)) \geq
\norm({\bf c}^+)$.
\end{myrules}
\end{definition}

\begin{lemma}\label{3.5}
$K$ has the halving property. 
\end{lemma}

\proof
We set $k(\halv({\bf c}^+)) = k'( 
\norm^{\frac{1}{2}}({\bf c}), k({\bf c}^+))\geq k({\bf c}^+)$ as in \ref{1.7}(4). 
Then we have that
$\norm(\halv({\bf c}^+))= f(\norm^{\frac{1}{2}}({\bf c}),k(\halv({\bf c}^+)))\geq
\norm({\bf c}^+)-1$, by  Definition~\ref{1.7}(4).

\nothing{
{\sf Why is $\sqrt{k({\bf c}^+) 
\cdot \norm^0({\bf c})} \geq k({\bf c}^+)$?
This is needed such that the halving gives stronger conditions.
}}

If ${\bf c}'$ is a simple creature and $\norm({\bf c}', k(\halv({\bf c}^+)) > 0$
and $\norm^{\frac{1}{2}}({\bf c}') \leq \norm^{\frac{1}{2}}({\bf c})$, 
then 
\begin{eqnarray*}
\norm({\bf c}',k({\bf c}^+)) &= &
f(\norm^{\frac{1}{2}}({\bf c}'),k({\bf c}^+))\\
\nothing{&=&
f(\norm^{\frac{1}{2}}({\bf c}),k(\halv({\bf c}^+))
\cdot \frac{k(\halv({\bf c}^+)}{k({\bf c}^+)})\\
&\geq &
f(\frac{\norm^{\frac{1}{2}}({\bf c}')}{k(\halv({\bf c}^+))})+  
f(\frac{k(\halv({\bf c}^+)}{k({\bf c}^+)})\\}
&\geq & f(\norm^{\frac{1}{2}}({\bf c}'),k(\halv({\bf c}^+)) +
f(\norm^{\frac{1}{2}}({\bf c}),k(\halv({\bf c}^+)))\\
&\geq & 1 + \norm({\bf c}^+) -1 \geq 
\norm({\bf c}^+).
\end{eqnarray*}
If  $\norm^{\frac{1}{2}}({\bf c}') > \norm^{\frac{1}{2}}({\bf c})$,
then the inequality follows from the monotonicity properties in Definition~\ref{1.7}
(4).
\proofend

\nothing{
{\sf This is stronger than property (3) of \ref{3.3}. 
But I think ``$-1$'' does not do any harm 
 as long as we find normal conditions and 
a fusion sequence in the proof of  \ref{3.8}.}
}

\begin{claim}\label{3.6}
Assume that $\name{\tau}$ is a $Q$-name for an ordinal, 
and let $a$ be a set of
ordinals. Let $m \in \omega$. Let $p$, $p^+$
be 
 conditions such that
\begin{myrules}
\item[(a)] $\eta = \rt(p)$,
\item[(b)] $(\forall \nu \in \dom(p)) (\norm({\bf c}_{p,\nu}, p(\nu)) \geq m)$,
\item[(c)]  $\dom(p^+) = \dom(p)$, and for $\eta \in \dom( p^+)$,
$k({\bf c}^+_{p^+,\eta})= k( \halv({\bf c}^+_{p,\eta}))$, 
${\bf c}_{p^+,\eta}=
{\bf c}_{p,\eta}$.
\end{myrules}
Then for any $q \in Q$:
$p^+ \leq q$ and $q \Vdash \name{\tau} \in a$ and $\nu=\rt(q)$ and
$\eta = \pr_{q,p}(\nu)$ imply that  there is some $q'$ such that
\begin{myrules}
\item[$(\alpha)$] $p^{\langle \eta \rangle} \leq q'$, $\nu = \rt(q')$
\item[$(\beta)$] $q' \Vdash \name{\tau} \in a$,
\item[$(\gamma)$] for every $\rho \in \dom(q')$, $\norm({\bf c}_{q',\rho},
q'(\rho))\geq m$.
\end{myrules}
\end{claim}

\nothing{{\sf What are $m^+$ and $a$?
I did not yet find this in [470]. In Chapter 2.2. there, halving 
and 2-bigness together with  some additional properties leads immediately
to properness. But there is is shown only for (not tree)-creatures.
Please, do not just cite. I do
not see how you can replace the 
``counting the finitely many possibilities for an 
extension of the stem''-arguments, because there are $\aleph_1$ partial
specialization functions allowed in the extension of
a condition (Which fall fortunately only into $\aleph_0$
isomorphism types. See below.) So far, \ref{3.6} is not used.}}
\proof 
\nothing{We define $p^+$ by $\dom(p^+) = \dom(p)$ and for
$\eta \in \dom(p)$, ${\bf c}^+_{p^+,\eta} =
\halv({\bf c}^+_{p,\eta})$. We prove that $p^+$ is as desired.}
So let $q \geq p^+$ and $q \Vdash \name{\tau} \in a$ and $\eta = 
\pr_{q,p}(\rt(q))$.
We take some $n(\ast) \in \omega$ such that
$$(\forall \rho \in \dom(q)) (\rho \in \bigcup_{n' \geq n(\ast)} q^{[n']}
\rightarrow \norm({\bf c}^+_{q,\rho}) > m).$$
We define $q'$ by $\dom(q') = \dom(q)$ and 
$\rho \in \bigcup_{n' \geq n(\ast)} q^{[n']}
\rightarrow {\bf c}^+_{q',\rho}={\bf c}^+_{q,\rho}$,
$\rho \in \bigcup_{n' < n(\ast)} q^{[n']}
\rightarrow {\bf c}^+_{q',\rho}=({\bf c}_{q,\rho},k({\bf c}^+_{p,\rho}))$.
\smallskip

$q$ and $q'$ force the same things, because we weakened $q$ to $q'$ only in 
an atomic part, because there are only finitely many $k$ such that
$({\bf c}_{q,\rho},k)$ is a creature with $0 \leq \norm^0({\bf c}_{q,\rho})$.
\smallskip

\relax From Lemma~\ref{3.5} we get $\rho \in \dom(q) \rightarrow \norm({\bf c}^+_{q',\rho}) 
\geq m$.
\proofend
\nothing{{\sf Starting from here until \ref{3.10}, 
things have to be checked very carefully
because Heike tried things on her own....}
\smallskip
\begin{definition}\label{3.7}
Let $\eta \in \spec$ and $p \in Q$ be smooth. 
Let the condition $(p,\eta)$ be undefined
if the norms of the following object are too small 
to make a condition.
Otherwise we  define $(p,\eta) = q$ by
$$\dom(q) = \{ \nu \cup  \eta \such 
(\forall \xi \in p) (\nu \subseteq \xi
\rightarrow \eta \cup \xi \in \spec)\}.$$
$q(\nu \cup \eta)  = p(\nu)$ for $\nu 
\cup \eta \in \dom(q)$, and
$\pr_{q,p}(\eta \cup \nu) = \nu \in p $.
\end{definition}
}

As a preparation for the following proof, we define isomorphism types of
partial specialization functions over conditions $p$:

\begin{definition}\label{3.7}
Let $\eta_0, \eta_1 \in \spec$ and let $p \in Q$.
We say $\eta_0$ is isomorphic  to $\eta_1$ over $p$ if there
is some injective partial function $f \colon T_A \to T_A$
such that $x <_{T_A} y$ iff $f(x) <_{T_A} f(y)$ and
$\dom(\eta_0) \cup \bigcup \{ \dom(\eta) \such  \eta \in \dom(p)\} \subseteq \dom(f)$
and $f \restriction \bigcup \{ \dom(\eta) \such  \eta \in \dom(p)\}
= id$ and $f[\dom(\eta_0)] = \dom(\eta_1)$ and 
$\eta_0(x) = \eta_1(f(x))$ for all $x \in \dom(\eta_0)$.
\end{definition}

Facts: For each fixed $p$, there are 
only countably many isomorphism types for $\eta$ 
over $p$. If the elements of $\dom(\eta_0)$ and of $\dom(\eta_1)$ are
pairwise incomparable in $T_A$ and if they
are isomorphic over $p$ with $\dom(p) = (T_A)_{<\alpha}$ for some countable $\alpha$,
and if there is some $r \geq p$ such that
$\eta_0 \in r^{[]}$, then there is some $r' \geq p$ such that
$\eta_1 \in (r')^{[]}$. 

\nothing{, and $(p,\eta_i)$ is a condition,
then $(p,\eta_0 \cup \eta_1)$ is a condition.
\begin{lemma}\label{3.8} (Some pure decision property)
Let $p \in P$ be smooth and let $\phi$ be a statement in the forcing
language. Then there is some $q \geq_0 p$ such that for all $r \geq q$,
if $r$ decides $\phi$, then there is some $\eta \in \dom(r)$ such that
$(q,\eta)$ decides $\phi$.
\end{lemma}
\proof 
We enumerate the isomorphism types of elements of $\spec$ over $p$ as
$\langle \eta_i \such i \in \omega \rangle$.
Now we build a sequence $\langle q_i \such i \in \omega \rangle$ such that
$q_{i+1} \geq_{i} q_i$ and all $q_i$ are normal.
We take $q_0 = p$.
Suppose $q_i$ has been chosen. If there is some normal
$r \geq q_i$ such that
$\rt(r)= \eta_i$ and $r$ decides the value of $\name{\tau}$,
we choose $\eta^\ast \in \dom(r)$ such that $\nor({\bf c}_{r,\eta^\ast} 
\geq i$.
We choose $\ell^*$ high enough and points $\xi \in \dom(q_i)$
 such that
for 
we can apply \ref{1.11} to ${\bf c}_{q_i,\xi}$ and 
$\eta^\ast$ for $\xi $ high enough.
Then
we graft $\eta $ onto $q_i$ in a way such that
the outcome $q_{i+1} \geq_{i} q_i$ and such that
$(q_{i+1}, \eta)$ has the same forcing power as $r$.
So we put in $\eta^\ast$ and above the gluing point we take $r$ instead of
$q_{i+1}$.
\proofend
\begin{lemma}\label{3.9a}
$Q$ is proper.
\end{lemma}
\proof Let $\tau$ be a $Q$-name for an ordinal and $p \in Q$.
We show how to find a countable set into which $\tau$ is forced.
First we take $q \geq p$ such that $q$ has the pure decision property.
Now we enumerate by $\langle \eta_i \such i < \omega_1 \rangle$ all
the finite partial specialization functions of $T_A$, such that
$(p,\eta_i)$ is a condition.
By the pure decision property we know that
in $V[G]$, $\tau[G] \in \{ \alpha \such 
(\exists i) (p,\eta_i) \Vdash \name{\tau} = \alpha \}$.
Assume that the latter set in not countable. 
Then there are $\aleph_1$ different $\alpha_i$ in it and 
hence $\aleph_1$ pairwise incompatible $(p,\eta_i)$.
We take the delta-lemma 
\cite[2, III, 5.4]{Jech} on the $\eta_i$ and get a root $\Delta$
and
$\aleph_1$  $\eta_i$ that all elements $\dom(\eta_i)
\setminus \dom(\Delta)$ are pairwise $<_{T_A}$ 
incomparable. Since there are only $\aleph_0$ isomorphy classes,
two of them are isomorphic over $p$. But then there is a 
condition stronger than $(p,\eta_i)$ and $(p,\eta_j)$, in contradiction to
that they force different values $\alpha_i$ to $\name{\tau}$.
\proofend
}

\begin{claim}\label{3.8}
Suppose that $p_0 \in Q$ and that $m < \omega$ and that $\name{\tau}$ 
is a $Q$-name of an ordinal. Then there is some $q \in Q$ such that
\begin{myrules}
\item[(a)]
$p_0 \leq_m q$,
\item[(b)] for some $\ell \in \omega$ we have that for every
$\eta\in q^{[\ell]} $ the condition $q^{\langle \eta\rangle}$
forces a value to $\name{\tau}$.
\end{myrules}
\end{claim}

\proof 
Choose $n(\ast)$ such that $\rho \in \bigcup_{n \geq n(\ast)}
p_0^{[n]} \rightarrow \norm({\bf c}^+_{p_0,\rho}) \geq m+1$.
Then we define $p_1$ by $\dom(p_1) = \dom(p_0)$ and 
$\rho \in \bigcup_{n' \geq n(\ast)} p_0^{[n']}
\rightarrow {\bf c}^+_{p_1,\rho}={\bf c}^+_{p_0,\rho}$,
$\rho \in \bigcup_{n' < n(\ast)} p_0^{[n']}
\rightarrow {\bf c}^+_{p_1,\rho}=\halv({\bf c}^+_{p_0,\rho})$.

\smallskip

Then we define 
\begin{equation*}
\begin{split}
X = \Biggl\{ \rho \such & \rho \in \bigcup_{n \geq n(\ast)} p_1^{[n]} \: 
\wedge (\exists q ) \biggl( p_1^{\langle \rho \rangle } \leq q 
\; \wedge \; q \mbox{ forces a value to } \name{\tau}\\
& \wedge (\forall \nu \in q^{[]}) (\norm({\bf c}^+_{q,\nu}) \geq 1)\biggr)\Biggr\}.
\end{split}
\end{equation*}

Let $p_2$ be chosen as in \ref{3.3} for $(p_1,X,n(\ast))$.
By a density argument, there is a front $\{\nu_0, \dots \nu_r\}$ of
$p_1$ such that for all $\nu_i$ the first clause of the alternative in \ref{3.3}(a) holds.

For $\tilde{m}<\omega$, $r,s \in Q$, $ \eta\in \dom(r)$, we denote 
the following property by $(\ast)^{\tilde{m},\eta}_{r,s}$:

\begin{equation}\tag*{$(\ast)^{\tilde{m},\eta}_{r,s}$}
\begin{split}&
r^{\langle \eta\rangle} \leq_0 s \; \wedge \\
&\forall \nu (\eta \subseteq \nu \in s \rightarrow 
\norm({\bf c}^+_{s,\nu}) \geq \tilde{m}+1) \: \wedge\\
&(\exists \ell \in \omega)(\forall \rho \in s^{[\ell]} ) (s^{\langle 
\rho \rangle} \mbox{ forces a value to } \name{\tau}).
\end{split}
\end{equation}

\smallskip

\nothing{Let $p_2 \leq_m q_0$ where $q_0$ is smooth. Such a $q_0$ can be chosen by
Lemma \ref{2.11}.}

We choose by induction on $t<\omega$ countable $N_t \prec ({\mathcal H}(\chi),\in)$ and
an ordinal $ \alpha_t$ and
pairs $(k_t,q_t)$ such that
\begin{myrules}
\item[(1)] $p_2, T_A, \name{\tau} \in N_0$, 
\item[(2)] $N_t \in N_{t+1}$,
\item[(3)] $N_t \cap \omega_1 = \alpha_t$,
\item[(4)] $\delta = \lim_{t \to \omega} \alpha_t$,

\item[(5)] $k_t$ is increasing with $t$, $k_t \geq n(\ast)$,
\item[(6)] $q_t\in Q$ is smooth, 

\item[(7)] $\alpha(q_{t}) = \alpha_t$,
\item[(8)] $k_t$ is the first $k$ strictly
larger than all the $k_{t_1}$ for $t_1 < t$
and such that $\rho \in q_t^{[\geq k]} 
\rightarrow \norm({\bf c}^+_{q_t,\rho}) > m +t+1$,
\nothing{and $(\forall n \leq t)(\forall \rho \in q_t^{[\geq k_]}) 
(e_n \in \dom(\rho))$ {\sf I added this in order to get
$\bigcup\{\dom(\eta) \such \eta \in r \} \supseteq T_{<\delta}$.
I do not see why the fusion of smooth conditions should be smooth.}
}
\item[(9)] $q_t \leq_{m+t+1} q_{t+1}$.

\item[(10)]if $\eta \in q_t^{[k_t]}$ and there is $q$ ($\in V$) satisfying
$(\ast)^{m+t+1,\eta}_{q_t,q}$, then
$q = q_{t+1}^{\langle \eta \rangle }$
satisfies it,

\item[(11)] $q_t \in N_{t+1}$,

\item[(12)] if $\eta \in q_t^{[k_t]}$ and no 
$q$ satisfies $(\ast)^{m+t,\eta}_{q_t,q}$,
then $(q_{t+1}^{\langle \eta\rangle}) ^{[]} = 
(q_{t}^{\langle \eta\rangle}) ^{[]}$ and $\eta \subseteq \rho \in q_t$
implies that
${\bf c}^+_{q_{t+1},\rho} = \halv({\bf c}^+_{q_t,\rho})$.
\end{myrules}

\nothing{
(13A)
An alternative (using just the simple creature part, and not
using halving) is replacing clause (viii) by
\begin{myrules}
\item[(viii')] if $\eta \in q_t^{[k_t]}$ and no 
$q$ satisfies $(\ast)^{m+t}_{q_t,\eta,q}$,
then  for every $\nu \in
(q_{t+t}^{\langle \eta\rangle}) ^{[]}$  there is no $r$ such that
 $(q_{t+1}^{\langle \nu\rangle})\leq_{m+t+1}
r$ and for some $\ell$ for every $\rho \in 
r^{[\ell]}$, 
$r^{\langle \rho\rangle }$ forces a value to $\name{\tau}$.
\end{myrules}
{\sf Saharon, what do you prove with this alternative?}
}

It is clear that the definition can be carried out 
as required. If we are given $q_t$ we can easily find $k_t$.
For each $\eta \in q_t^{k_t}$ we choose $q_{t,\eta} \in N_t$ such that
$(\ast)^{m+t+1}_{q_t,q_{t,\eta}}$ if possible 
and in fact w.l.o.g. $q_{t,\eta} = q_t^{\langle \eta \rangle}$,
otherwise we follow (12) and apply the halving function.

\smallskip

Having carried out the induction, we let $r = 
\bigcup_{t \in \omega} (q_t\restriction q_t^{(k_{t-1},k_t]})$.
So, by (7), 
$r \in Q$ is smooth with $\alpha(r)= \delta$ and for every $t$ we have
$q_t \leq_{m+t+1} r$, and in particular $p_2 \leq_{m+1} r$.

\nothing{
Let
$$X= \{ \rho \in r^{[]} \such (\exists \ell) 
(\forall \nu \in (r^{\langle \rho \rangle})^{[\ell]})
(r^{\langle \nu \rangle} \mbox{ forces a value to } \name{\tau})\}.$$
Notice that $X$ is upwards closed, and so we can apply 
\ref{3.3} with $(r,X,m) $ for $(p,X,k^*)$ and get a condition $q_\ast$
as there.}

Assume for a contradiction that we are in the bad case 
\begin{equation}\label{otimes}\tag*{$\otimes$}
(\forall \ell \in \omega)(\exists \rho^\ast \in r^{[\ell]})
(r^{\langle \rho^\ast \rangle} \mbox{ does not decide the value of } \tau).
\end{equation}

Choose a minimal $\ell$ as in $\otimes$ and a $\rho^\ast$ as there.
Choose 
$q \geq r^{\langle \rho^\ast\rangle}$
in $Q$ such that $q$ forces a value to $\name{\tau}$.

\nothing{We should prove that $\rho \in X$. Then we have 
$\rho \in (p_2)^{\langle \rho^\ast\rangle}\cap X$ and $\otimes$ is 
impossible.}

Let $\nu_0 \in q^{[]}$ be such that 
$$ \nu_0 \subseteq \nu \in q \rightarrow \norm({\bf c}^+_{q,\nu}) \geq m +1.$$
W.l.o.g.\ (otherwise we strengthen $q$) we assume that $\nu_0=\rt(q)$.

\smallskip

As $r$ is smooth,
by the definition of $q \geq r$ (\ref{2.2}(f)) we have that 
the additional information on partial specialization functions
that are in $q$ but not is $r$ does not have the domain in
 $(T_A)_{<\delta}$.

\smallskip

\nothing{{\sf Starting from here, I just typed (including 13A, 14A), 
because I understand less and less.}}

Let $t(\ast)$ be such that $\dom(\rt(q)) \cap (T_A)_{<\delta} 
\subseteq T_{< \alpha_{t(\ast)}}$, and w.l.o.g.\
$\nu_0 \in q^{[k_{t(\ast)} + i(p_0) -i(q)]}$,
so $\pr_{q,r}(\nu_0) \in r^{[k_{t(\ast)}]}$, $i(p_0) =
i(q_t) = i(r)$.
Now easily 
\begin{equation}
\tag{$\ast$}
\begin{split}
&\mbox{ if } t(\ast) \leq t, \nu_0 \leq \nu \in \dom(q), 
\pr_{q,r}(\nu) \in r^{[k_t]}, \\
&\mbox{ then } \pr_{q,r}(\nu) = \pr_{q,q_t}(\nu), \dom(\nu) \cap
 (T_A)_{\alpha(r)} = \dom(\pr_{q,r}(\nu)) =\\
&
\dom(\pr_{q,q_t}(\nu)) \subseteq (T_A)_{< \alpha(q_t)}.
\end{split}
\end{equation}

\nothing{Now $\eta_0 = \pr_{q,r}(\nu_0)$ cannot be in $X$
as $\nu_0 \in q \subseteq q_\ast^{\langle \rho^* \rangle} $ and
$q_\ast^{\langle \rho^\ast \rangle} \cap X = \emptyset$ and 
$X$ is upwards closed and $\eta_0 \subseteq \nu_0$.}
\nothing{
(as $\pr_{q,r}(\nu_0)$ necessarily 
$\in \dom(q_0) \subseteq \dom(r) X$ by the assumption
towards a contradiction above).
}

Now $\eta_0 \in r^{[k_{t(\ast)}]}$ and even $\eta_0 \in 
q_{t(\ast)}^{[k_{t(\ast)}]}$. So by the choice of
$\langle q_t \such t \in \omega \rangle$ we know that there is no $q$
with $(\ast)^{m+t(\ast)+1,\eta_0}_{q_{t(\ast)},q}$,
as otherwise $q_{t(\ast)+1}^{\langle \eta_0\rangle}$ 
would be like this and this property
would be inherited by $r$.
So clause (12)  applies, which means 

\begin{equation}\label{box}
\tag{$\boxtimes$}
\begin{split}
& \mbox{ there is no $\rho$ such that } 
\eta_0  \subseteq \rho \in \dom(q_{t(\ast)})
\mbox{ and there is no } s \mbox{ such that } \\
& q_{t(\ast)}^{\langle \rho \rangle} \leq_{m+t(\ast)+1} s
\mbox{ and } (\exists \ell)(\forall \rho' \in s^{[\ell]})
(s^{\langle \rho'\rangle} \mbox{ forces a value to } \name{\tau}).
\end{split}
\end{equation}

Choose $\nu^\ast$ such that $\nu_0 \subseteq \nu^\ast \in\dom(q)$ and
$\pr_{q,r}(\nu^\ast) \in r^{[k_{t(\ast)}]}$.
Let $\eta^\ast = \pr_{q,q_t}(\nu^\ast) \in q_t^{[k_{t(\ast)}]}$.

\nothing{{\sf Why can we not just take $\nu^\ast = \nu_0$?}}

Fix for some time $\nu \in \suc_q(\nu^\ast)$ and let
$\eta = \pr_{q,q_t}(\nu) \in q_t^{[k_t]}$, so $\eta \in
\suc_{q_t}(\eta^\ast)$.

So let $\dom(\nu) \setminus \dom(\eta) = 
\{ x_0, \dots, x_{\tilde{s}-1}\}$. We just saw that
$x_0, \dots , x_{\tilde{s}-1} \not\in (T_A)_{<\delta}$.
Let us define
$\bar{y}= \langle y_\ell \such \ell < \tilde{s} \rangle$ is a candidate
for an extended domain iff:
\begin{myrules}
\item[(a)] $y_\ell$ are without repetitions,
\item[(b)] there is some $r_{\bar{y}}$ such that
\begin{myrules}
\item[(0)] $
\rt(r_{\bar{y}})
  = \eta \cup \{(y_\ell,\nu(x_\ell)) \such \ell < \tilde{s} \} \in Q$,
such that

\item[(1)] $r_{\bar{y}} \geq q_t^{\langle \eta \rangle}$,
\item[(2)] $(\forall \rho)(\eta \subseteq \rho \in \rt(r_{\bar{y}})
\rightarrow \norm({\bf c}^+_{\rho,r_{\bar{y}}}) > m+
t(\ast) +1)$,
\item[(3)] $r_{\bar{y}}$ forces a value to $\name{\tau}$
\item[(4)] $\langle (y_\ell,\nu(x_\ell)) \such \ell < \tilde{s} \rangle 
\in \spec^{T_A}$ is isomorphic over $(T_A)_{<\delta}$ to

$\langle (x_\ell,\nu(x_\ell)) \such \ell < \tilde{s} \rangle$.

\nothing{ for $\ell < s$, $y_\ell \restriction \alpha_{t(\ast)}
=x_\ell \restriction\alpha_{t(\ast)}$. {\sf Here, we think of 
$x_\ell$, $y_\ell \in \omega_1$ not as ordinals, but as 
branches in the Aronszajn tree, $x_\ell = \{ x \in T_A \such x 
<_{T_A} x_\ell\}$. Is this interpretation intended?}}
\end{myrules}
\end{myrules}
We set $$Y = Y_\eta =
\{ \bar{y} \such \bar{y} \mbox{ is a candidate for an extension}\}.$$

Now we have that $\langle x_\ell \such \ell < \tilde{s} \rangle \in Y$. 
This is exemplified
by $q^{\langle \eta \rangle }$.
\nothing{ which we assume to be
treated with the operation $^+$ as in Claim \ref{3.6}
for $m + t(\ast)+1$.}

Now we have that $q_t \in N$ and for all $\ell$, $x_\ell \in (T_A)_{\geq 
\delta}$, because the $\alpha_t$ are cofinal in $\delta$ and since
$\alpha_t = \alpha(q_t)$.

\smallskip
Since $x_\ell \geq \delta$, counting isomorphism types over $(T_A)_{<\delta}$
 yields $|Y_\eta| = \aleph_1$. 

\smallskip
 By a fact on Aronszajn trees (Jech,  or \cite[III, 5.4]{Sh:h}) we find
$\langle y^\eta_{j,\ell} \such \ell <s, j \in \omega_1
\rangle$ and a root $\Delta_\eta$ such that

\begin{myrules}
\item[(a)] $y^{\eta}_{j,\ell} \in Y_\eta$ are without repetition,
\item[(b)] for $j \neq j'$, $\{ y_{j,\ell}^{\eta} \such \ell < \tilde{s} \} \cap
\{  y_{j',\ell}^{\eta} \such \ell < \tilde{s} \} = \Delta_\eta$,
\item[(c)] if $j_1 \neq j_2$ and if  $y^\eta_{j_1,\ell_1} \not\in \Delta_\eta$
and $y^\eta_{j_2,\ell_2}\not\in \Delta_\eta$ then
they  are incompatible in $<_{T_A}$.
\end{myrules}

Let $r^\eta_{\langle y_{j,\ell}^{\eta} \such \ell <\tilde{s} \rangle}$
witness that $\langle y_{j,\ell}^{\eta} \such \ell <\tilde{s} \rangle \in Y_\eta$.

Let ${\bf c} = \{ \pr_{q,q_t}(\nu) \such \nu 
\in 
\suc_q(\nu^*)\}$.
This 
is a simple $(i(q_t) + k_{t(\ast)} +1)$-creature 
with $\norm({\bf c}, q_t(\eta^*))
\geq m + t(\ast) +2$ by property (10) of $(k_{t(\ast)},q_{t(\ast)})$.
For each $\eta \in \rge(\val ({\bf c})) $ let 
$\langle y_{j,\ell}^{\eta} 
\such  \ell < \tilde{s}, j< \omega_1 \rangle$ be
as above and 
let $r^\eta_j$ be a witness for 
$\langle y^\eta_{j,\ell} \such \ell < \tilde{s} \rangle \in Y_\eta$.

\smallskip

Let $j^\ast = \norm^0({\bf c}_{q_{t(\ast)},\eta^*})$.

\nothing{
Choose $y \in T_\delta \setminus \{ x_{j,\ell}^{\eta,*} \such \ell <s,
 j < j^*\}$.
For each $\eta \in \rge(\val({\bf c}))$ let $\bar{k_\eta}
= \langle k_{\eta,j} \such j<j^*\rangle$ be such that
\begin{itemize}
\item
$\bar{k_\eta}$ is without repetition,
\item $k_{\eta,j} < n_{2,i(q_t) + k_t -1}$
\item $\eta \cup \{(y,k_{\eta,j}) \} \in \spec$.
\end{itemize}
We can find such $\bar{k}$ since $|\dom(\eta)| 
<< n_{2,i(q_t) +k_t -1}$.}

For each $\eta \in \rge(\val({\bf c}_{q_{t(\ast)},\eta^\ast}))$ choose a witness 
$\nu_\eta \in \suc_q(\nu^*)$ such that $\pr_{q,q_{t(\ast)}}(\nu_\eta) 
= \eta$.
Now we define a simple $i({\bf c}_{q_{t(\ast)},\eta^*})$-creature ${\bf d}$ by

$$\eta({\bf d})= \eta({\bf c}_{q_{t(\ast)},\eta^*})$$
$$\rge(\val({\bf d}))= \{\eta \cup
\{(x_{j,\ell}^\eta,\nu_\eta(x_\ell)) \such \ell < \tilde{s} \}
\such \eta\in \rge(\val({\bf c})), j < j^*\}.$$

Then we have by Claim~\ref{1.9} that ${\bf d}$ is a $i({\bf c}_{q_{t(\ast)},
\eta^*}))$-creature and $\norm^0({\bf d}) = 
\norm^0({\bf c}_{q_{t(\ast)},\eta^\ast})=m+t(\ast) +1$.
$\norm^\ast$ drops at most by 1. So we have
$\norm({\bf d},q_{t(\ast)+1}(\eta^\ast)) >0$ and hence by Claim~\ref{3.6}
$\norm({\bf d},q_{t(\ast)}(\eta^\ast)) >m$.
Now we define $s \in Q$ as follows:
\begin{myrules}
 \item[$(\alpha)$] $\rt(s) = \eta^*$, $s(\eta^*) = q_{t(\ast)}(\eta^*)$,
\item[$(\beta)$] ${\bf c}_{s,\eta^*} = {\bf d}$,
\item[$(\gamma)$] if $\rho \in \rge(\val({\bf d}))$ and if
$\rho = \eta \cup
\{(x_{j,\ell}^\eta,\nu_\eta(x_\ell)) \such \ell < \tilde{s} \}
\}$
then $s^{\langle \rho\rangle} = r_j^\eta$.
\end{myrules}

\smallskip

Clearly $s \in Q$ and $q_t^{\langle \eta^*\rangle} \leq_{m+t(\ast) +1}
s$ 
\nothing{
{\sf for this property, $\norm$ is used and not $\norm^0$. The result with
$\norm^0$ seems to be even stronger. Where is halving used?
Where do you have to lessen to $\norm$ and work with non-simple creatures?}
}
and 
for every
$\eta\in q^{[\ell]} $ the condition $q^{\langle \eta\rangle}$
forces a value to $\name{\tau}$, in fact $\ell = k_{t(\ast)}+1$ 
is o.k., by the way the $r_j^\eta$ were chosen. So we get a contradiction to
$\boxtimes$
\nothing{({\sf 
which symbol I did not find. 
I think it is $q_\ast^{[]} \cap X= \emptyset$})}
 and to the choice of $\eta^*$.
\proofend

\def\germ{\frak} \def\scr{\cal}
  \ifx\documentclass\undefinedcs\def\rm{\fam0\tenrm}\fi
  \def\defaultdefine#1#2{\expandafter\ifx\csname#1\endcsname\relax
  \expandafter\def\csname#1\endcsname{#2}\fi} \defaultdefine{Bbb}{\bf}
  \defaultdefine{frak}{\bf} \defaultdefine{mathbb}{\bf}
  \defaultdefine{mathcal}{\cal}
  \defaultdefine{beth}{BETH}\defaultdefine{cal}{\bf} \def\bbfI{{\Bbb I}}
  \def\mbox{\hbox} \def\text{\hbox} \def\om{\omega} \def\Cal#1{{\bf #1}}
  \def\pcf{pcf} \defaultdefine{cf}{cf} \defaultdefine{reals}{{\Bbb R}}
  \defaultdefine{real}{{\Bbb R}} \def\restriction{{|}} \def\club{CLUB}
  \def\w{\omega} \def\exist{\exists} \def\se{{\germ se}} \def\bb{{\bf b}}
  \def\equivalence{\equiv} \let\lt< \let\gt> \def\cite#1{[#1]}


\end{document}